\documentclass[11pt]{article} 
\usepackage[a4paper]{geometry}
\usepackage{bookman,fourier}
\usepackage{amssymb,amsbsy,amsmath}

\usepackage{notestyle}
\usepackage{epstopdf}

\begin{document}
\title{Generalized Taylor operators and Hermite subdivision schemes}

\author{Jean-Louis Merrien\footnote{INSA de Rennes, IRMAR - UMR 6625,
    20 avenue des Buttes de Coesmes, CS 14315, 35043 RENNES CEDEX,
    France, email: \texttt{jmerrien@insa-rennes.fr}}
  \and
  Tomas Sauer\footnote{Lehrstuhl f\"ur Mathematik mit Schwerpunkt
    Digitale Bildverarbeitung \& FORWISS, Universit\"at Passau,
    Fraunhofer IIS Research Group ``Knowledge Based Image Processing'',
    Innstr. 43, D-94032 PASSAU, Germany, email:
    \texttt{Tomas.Sauer@uni-passau.de}}} 

\maketitle

\begin{abstract}
  Hermite subdivision schemes
  act on vector valued data that is
  not only considered as functions values in $\RR^r$, but as
  consecutive derivatives,
  which leads to a mild form of level dependence of the scheme. 
  Previously, we have proved that a
  property called spectral condition or sum rule implies a
  factorization in terms of a generalized difference operator that gives
  rise to a ``difference scheme'' whose contractivity governs the
  convergence of the scheme. 
  But many convergent Hermite schemes, for example, those based on
  cardinal splines, do not satisfy the spectral condition.
  In this paper, we generalize the property in a way that preserves
  all the above advantages: the associated
  factorizations and convergence theory. Based on these results, 
  we can include the case of cardinal splines and also enables us to
  construct new types of convergent Hermite subdivision schemes.
\end{abstract}

\section{Introduction}
Subdivision schemes, as established in \cite{CDM91}, are efficient
tools for building curves and 
surfaces with applications in design, creation of images and motion
control. For vector subdivision schemes,
cf. \cite{DynLevin02,Han03,MicchelliSauer98}, it is not so straightforward 
to prove more than the H\"older regularity of the limit function, due
to the more complex nature of the underlying factorizations. On the
other hand, Hermite subdivision schemes
\cite{dyn95,han03:_desig_hermit,HanXueYu05,jeong17:_const_hermit,ManniGuglielmi2011}
produce function vectors 
that consist of consecutive derivatives of a certain function, so that
the notion of
convergence automatically includes regularity of the leading component
of the limit. Such schemes have even been considered also on manifolds
recently \cite{moosmueller16:_c_analy_hermit_subdiv_schem_manif} and
have also been used for wavelet constructions
\cite{cotronei17:_hermit}. While vector subdivision schemes are
quite well--understood, nevertheless there are still surprisingly many open
questions left in Hermite subdivision. In particular, a
characterization of convergence in terms of factorization and
contractivity is still missing.

In previous papers
\cite{DubucMerrien09,MerrienSauer2012,MerrienSauer2017}, we
established an equivalence between a so--called \emph{spectral
  condition} and operator 
factorizations that transform a Hermite scheme into a vector scheme
for which analysis tools are available. Under this transformation, the
usual convergence of the vector subdivision scheme implies convergence
for the Hermite scheme and thus regularity of the limit function. It
was even conjectured for some time that the spectral condition,
sometimes also called the \emph{sum rules}
\cite{conti16:_approx,HanXueYu05} of the Hermite subdivision
scheme, might be necessary for convergence.

In this paper we show, among others results, that this conjecture does not
hold true. We
define a new set of significantly more general spectral conditions,
called \emph{spectral chains}, 
that widely generalize the classical spectral condition from
\cite{DubucMerrien09} and show that these spectral conditions are more
or less equivalent to the existence of a factorization with respect to
respective generalized Taylor operators and allow for a description of
convergence by means of contractivity. We then define a process that
allows us to construct Hermite subdivision schemes of arbitrary
regularity with guaranteed convergence and, in particular, give
examples of convergent Hermite subdivision schemes that do not satisfy
the spectral condition. In addition, 
our new method can be
applied to an example based on B--splines and their derivatives which
was one of the first examples of a convergent Hermite subdivision
scheme that does not satisfy the spectral condition,
\cite{MerrienSauer2011}.

The paper is organized as follows: after introducing some basic
notation and the concept of convergent vector and Hermite subdivision
schemes, we introduce the new concept of chains and generalized Taylor
operators in Section~\ref{D:SpecVecs} and 
use them 
for the factorization of subdivision operators in
Section~\ref{sec:chainsfakt}. These results allow us to extend the
known results about the convergence of the Hermite subdivision schemes
to this more general case in Section~\ref{sec:Convergence}. 
Section~\ref{sec:unfactor} 
is devoted to the construction of 
a convergent Hermite subdivision scheme emerging from a
properly constructed
contractive vector subdivision scheme by reversing the factorization
process, even in the generality provided by generalized Taylor
operators.
Finally, we give some examples of the results of such constructions in
Section~\ref{sec:Examples}, and also provide a new approach for the
aforementioned spline case.

\section{Notation and fundamental concepts}
Vectors in $\RR^r$, $r \in \NN$, will 
generally 
be labeled by lowercase boldface letters: $\yb =
\left[ y_j \right]_{j=0,\dots,r-1}$
or  $\yb =\left[ y^{(j)} \right]_{j=0,\dots,r-1}$, where 
the
latter notation is used to highlight the fact that in Hermite subdivision the
components of the vectors correspond to derivatives. 
Matrices in $\RR^{r  \times r}$
will be written as uppercase boldface letters,
such as $\Ab = \left[ a_{jk} \right]_{j,k=0,\dots,r-1}$.
The space of polynomials in one variable of degree at most $n$ will be
written as $\Pi_n$, while $\Pi$ will denote the space of all polynomials. 
Vector sequences will be considered as
functions from $\ZZ$ to $\RR^r$ and the vector space of all such
functions will be denoted by $\ell(\ZZ,\RR^r)$ or $\ell^r ( \ZZ)$.
For $\yb(\cdot)\in\ell(\ZZ,\RR^r)$, the \emph{forward difference} is
defined as  $\Delta\yb(\alpha):=\yb(\alpha+1)-\yb(\alpha)$, $\alpha
\in \ZZ$, and iterated to $\Delta^{i+1}\yb := \Delta \left( \Delta^i
  \yb \right) =
\Delta^{i}\yb(\cdot+1)-\Delta^{i}\yb(\cdot)$,  $i\ge 0$.

Given a finitely supported sequence of matrices $\Ab = \left( \Ab
  (\alpha) \right)_{\alpha \in \ZZ} \in \ell^{r \times r} (\ZZ)$,
called the \emph{mask} of the subdivision scheme, we define the
associated \emph{stationary subdivision operator}
$$
S_\Ab : \cb \mapsto \sum_{\beta \in \ZZ} \Ab ( \cdot - 2\beta ) \, \cb
(\beta), \qquad \cb \in \ell^r (\ZZ).
$$
The iteration of subdivision operators $S_{\Ab_n}$, $ n\in\NN$, is
called a \emph{subdivision scheme} and consists of the successive
applications
of level-dependent subdivision operators, acting on vector valued data,
$S_{\Ab_n} : \ell^r \left( \ZZ \right) \to \ell^r\left( \ZZ \right)$,
defined as 
\begin{equation}
  \label{eq:SubdOpDef}
  \cb_{n+1}(\alpha) =S_{\Ab_n} \cb_{n}(\alpha) := \sum_{\beta \in \ZZ}
  \Ab_n \left( \alpha - 2 
    \beta \right) \, \cb_n (\beta), \qquad \alpha \in \ZZ, \qquad \cb
  \in \ell^r \left( \ZZ\right).
\end{equation}
An important algebraic tool for stationary subdivision operators is
the \emph{symbol} of the mask,
which is 
the matrix valued Laurent polynomial
\begin{equation}
  \label{def:Laurentpoly}
  \Ab^* (z) := \sum_{\alpha \in \ZZ} \Ab (\alpha) \, z^\alpha, \qquad z
  \in \CC \setminus \{ 0 \}.
\end{equation}
We will focus our interest on two kinds of such schemes, the
first one being ``traditional'' vector subdivision schemes in the sense of
\cite{CDM91}, where $\Ab_n$ is independent of $n$, i.e.,
$\Ab_n(\alpha)=\Ab(\alpha)$ for any $\alpha\in\ZZ$ and any $n\ge 0$.
In the following, such schemes for which an elaborate theory of
convergence exists, will simply be called a {\em vector scheme}. Their
convergence is defined in the following way.

\begin{definition}\label{def:convergvector}
Let $S_\Ab :
\ell^r ( \ZZ )\to \ell^r ( \ZZ )$
be a vector subdivision operator. The operator is 
\emph{$C^p$--convergent}, $p \ge 0$, if for any data $\gb_0\in\ell^r ( \ZZ)$ 
and corresponding sequence of refinements 
$\gb_n = S_\Ab^n \gb_0$ there exists a function $\psi_\gb \in C^p \left(
  \RR,\RR^r \right)$ 
such that for any compact $K\subset\RR$ there exists a sequence $\varepsilon_n$ 
with limit $0$ that satisfies 
\begin{equation}
  \label{eq:StatC0Conv}
  \max_{\alpha \in \ZZ \cap 2^n K} \left\| \gb_n (\alpha) - \psi_\gb \left( 2^{-n}
      \alpha \right) \right\|_\infty \le \varepsilon_n.  
\end{equation}
\end{definition}

\noindent
As the second type of, now even level--dependent, schemes we consider
the {\em Hermite scheme} where
$\Ab_n(\alpha)=\Db^{-n-1}\Ab(\alpha)\Db^n$ for $\alpha\in\ZZ$ and
$n\ge 0$
with 
the diagonal matrix 
$\Db :=  \begin{bmatrix}
  1 \\
  & \frac12 \\
  & & \ddots \\
  & & & \frac1{2^d} \
\end{bmatrix}$. In this case $r=d+1$ and for $k=0,\ldots,d$
the k-th component of $\cb_n(\alpha)$ corresponds to an approximation
of the  k-th derivative of some function 
$\varphi_n$ 
at $\alpha 2^{-n}$.
Starting from an initial sequence $\cb_0$, a Hermite scheme can be rewritten
\begin{equation}
  \label{eq:defHermop}
  \Db^{n+1} \, \cb_{n+1}(\alpha)
  = \Db^{n+1} \,S_{\Ab} \Db^n \cb_n (\alpha)
  = \sum_{\beta\in\ZZ} \Ab
  \left(\alpha - 2 \beta \right) \, \Db^n  \, \cb_n(\beta), \qquad
  \alpha \in \ZZ,\quad n\ge 0.
\end{equation}
Convergence of Hermite schemes is a little bit more intricate and
defined as follows.

\begin{definition}\label{def:HermCdConv}
  Let $\Ab \in \ell^{(d+1)\times(d+1)} (\ZZ)$ be a mask and
  $H_\Ab$ the associated Hermite subdivision scheme
  on $\ell^{d+1}(\ZZ)$
  as defined in \eqref{eq:defHermop}. The scheme
  is \emph{convergent} if for
  any data $\fb_0 \in \ell^{d+1} ( \ZZ)$ and the corresponding
  sequence of refinements $\fb_n=[f_n^{(0)} ,\ldots,f_n^{(d)}]^T$,
  there exists a function 
  $\Phi=[\phi_{i}]_{0\le i\le d}\in C\left(\RR,\RR^{d+1} \right)$
  such that for any compact $K\subset\RR$ there exists a sequence
  $\varepsilon_n$ with limit $0$ which satisfies 
  \begin{equation}
    \label{eq:HermCdConv}
    \max_{0\le i\le d}\max_{\alpha \in \ZZ \cap 2^n K} \left|
      f_n^{(i)} (\alpha) -\phi_{i} \left( 2^{-n} \alpha \right)
    \right| \le \varepsilon_n.
  \end{equation}
  The scheme
  $H_\Ab$ is said to be \emph{$C^{p}$--convergent} with $p\ge d$ if moreover
  $\phi_0 \in C^{p}\left( \RR, \RR \right)$ and
  $$
  \phi_0^{(i)}=\phi_i,\quad 0\le i\le d.
  $$ 
\end{definition}

\begin{remark}
  Since the intuition of Hermite subdivision schemes is to iterate on
  function values and derivatives, it usually only makes sense to
  consider $C^p$--convergence for $p \ge d$. Note, however, that the
  case $p > d$ leads to additional requirements.
\end{remark}

\section{Generalized Taylor operators and chains}
\label{sec:GenTayl}
In this section, we introduce the concept of generalized Taylor
operators and show that they form the basis of symbol
factorizations. The first definition concerns vectors of almost monic
polynomials of increasing degree.

\begin{definition}\label{D:SpecVecs}
  By $\VV_d$ we denote the set of all vectors $\vb$ of polynomials in
  $\Pi_{d}$ with the property that
  \begin{equation}
    \label{eq:SpecVecs}
    \vb = \left[
      \begin{array}{c}
        v_d \\ \vdots \\ v_0
      \end{array}
    \right], \qquad v_j = \frac{1}{j!} (\cdot)^j + u_j \in \Pi_j,
    \qquad u_j \in \Pi_{j-1}.
  \end{equation}
  A vector in $\VV_d$ thus consists of polynomials $v_j$ of degree
  \emph{exactly} $j$ whose leading coefficient is normalized to
  $\frac1{j!}$, and the remaining part of the polynomial $v_j$ of
  lower degree is denoted by $u_j$.
\end{definition}

\noindent
Note that in \eqref{eq:SpecVecs} we always have $v_0 = 1$ and $u_0 = 0$. Also keep
in mind that the vectors $\vb$ are {\em indexed in a reversed order}, but
referring directly to the degree of the object, this notion is more
comprehensible.

We will use the convenient notation of \emph{Pochhammer symbols}
$(\cdot)_j \in \Pi_j$, $j \ge 0$, in the following way:
\begin{equation}
  \label{eq:Pochhammer}
 (\cdot)_0 :=1, \qquad (\cdot)_j := \prod_{k=0}^{j-1} (\cdot
 - k ),\quad j\ge 1, \qquad \text{ and } \qquad
 [\cdot]_j := \frac{1}{j!} (\cdot)_j,\quad j\ge 0.
\end{equation}
These polynomials satisfy
\begin{equation}
  \label{eq:PochhammerProps}
  \Delta (\cdot)_j = j \, (\cdot)_{j-1}, \qquad
  \Delta [\cdot]_j = [\cdot]_{j-1}.
\end{equation}
Both $\left\{ (\cdot)_0,\dots,(\cdot)_j \right\}$ and
$\left\{ [\cdot]_0,\dots,[\cdot]_j \right\}$ are bases of $\Pi_j$ and
allow us to write the
Newton interpolation formula of degree $d$ at $0,\dots,d$ in the form
$$
x^j = \sum_{k=0}^j \frac{1}{k!} \left( \Delta^k (\cdot)^j \right) (0) \,
(x)_k = \sum_{k=0}^j \left( \Delta^k (\cdot)^j \right) (0) \, [x]_k;
$$
then, since $\Delta^j (\cdot)^j = j!$, we have that
$$
\frac1{j!} (\cdot)^j = [\cdot]_j + \sum_{k=0}^{j-1} \frac{\left(
    \Delta^k (\cdot)^j \right) (0)}{j!} \, [x]_k
$$
which implies that
\begin{equation}
  \label{eq:Pochv}
  \vb \in \VV_d \qquad \Leftrightarrow \qquad v_j = [\cdot]_j + u_j,
  \qquad u_j \in \Pi_{j-1}\quad j=0,\ldots,d.
\end{equation}
We will use this form in the future to write each $\vb \in \VV_d$ as
\begin{equation}
  \label{eq:vwFormula}
  \vb = \left[
    \begin{array}{c}
      {[\cdot]_d} \\ \vdots \\ {[\cdot]_0}
    \end{array}
  \right] + \ub.  
\end{equation}
Generalizing the Taylor operators operating on vector functions
$\RR\to\RR^{d+1}$ introduced in
\cite{DubucMerrien09,MerrienSauer2012}, we define the following
concept.

\begin{definition}
  A \emph{generalized incomplete Taylor operator} is an operator of
  the form
  \begin{equation}
    \label{eq:TaylorGen}
    T_d :=\begin{bmatrix} 
      \Delta & -1 & * & \dots & * \\
      & \ddots & \ddots & \ddots & \vdots \\
      &  & \ddots & \ddots & * \\
      & & & \Delta & -1 \\
      & & & & 1
    \end{bmatrix}
    = \begin{bmatrix} 
      \Delta I & \\ & 1
    \end{bmatrix}  + \left[ t_{jk} \right]_{j,k = 0,\dots,d},
  \end{equation}
  where $t_{j,j+1} = -1$ and $t_{jk} = 0$ for $k \le j$. In the same
  way, the \emph{generalized complete Taylor operator} is of the form
  \begin{equation}
    \label{eq:TaylorGenComp}
    \widetilde T_d := \begin{bmatrix} 
      \Delta & -1 & * & \dots & * \\
      & \ddots & \ddots & \ddots & \vdots \\
      &  & \ddots & \ddots & * \\
      & & & \Delta & -1 \\
      & & & & \Delta
    \end{bmatrix}  = \Delta I + \left[ t_{jk}\right]_{j,k = 0,\dots,d}.
  \end{equation}
\end{definition}

\begin{remark}
 The Taylor operator becomes generalized for $d \ge 2$, otherwise we simply
 recover the classical case, see Example \ref{Ex:ChainEx}. 
\end{remark} 

\begin{lemma}\label{L:TaylorKernels}
  Let $\vb:=[v_d,\ldots,v_0]^T$ be a vector of polynomials in
  $\Pi^{d+1}$ with $v_0 = 1$. Then $\vb \in \VV_d$ if and only if
  there exists a generalized complete Taylor operator $\widetilde T_d$
  such that $\widetilde T_d \vb = 0$.
\end{lemma}

\begin{pf}
  For ``$\Leftarrow$'' suppose that $\widetilde T_d \vb = 0$ and let
  us prove inductively
  for   $j=0,\dots,d$, that
  $v_j = [\cdot]_j + u_j$, for appropriate $u_j \in \Pi_{j-1}$.
  The assumption $v_0 = 1$ ensures that for $j=0$ by simply setting
  $u_0=0$.
  Now, for $0 \le j < d$, we assume that $v_{j+1}$ is of degree $m \ge
  0$ and write it in the basis $\left\{ [\cdot]_0,\ldots, [\cdot]_m
  \right\}$ as 
  $$
  v_{j+1} = \sum_{k=0}^m c_k [\cdot]_k = \sum_{k=j+2}^m
  c_k [\cdot]_k + c_{j+1} [\cdot]_{j+1} + q,
  $$
  with $q \in \Pi_j$, hence $\Delta q \in \Pi_{j-1}$.
  By induction, we suppose that $v_j = [\cdot]_j + u_j$, $u_j \in
  \Pi_{j-1}$ and $v_k\in\Pi_k$ for $k=0,\ldots,j-1$. Then $\widetilde
  T_d \vb = 0$ implies at row $d-j-1$ that
  \begin{eqnarray*}
    0 & = & \Delta v_{j+1} - v_j + \sum_{k=0}^{j-1} t_{d-j-1,d-k} v_k \\
      & = & \sum_{k=j+2}^m c_k [\cdot]_{k-1} + c_{j+1}
            [\cdot]_j + \Delta q - [\cdot]_j - u_j +
            \sum_{k=0}^{j-1} t_{d-j-1,d-k}v_k\\
      & = & \sum_{k=j+1}^{m-1} c_{k+1} [\cdot]_k + \left(
            c_{j+1} - 1 \right) [\cdot]_j + u, \qquad u\in \Pi_{j-1},
  \end{eqnarray*}
  and comparison of coefficients yields $c_{j+2} = \cdots = c_m = 0$
  as well as $c_{j+1} = 1$, hence $v_{j+1} = [\cdot]_{j+1} + u_{j+1}$
  with $u_{j+1} \in \Pi_j$, which advances the induction hypothesis.

  For the converse ``$\Rightarrow$'',  we note that for any $\vb \in
  \VV_d$ we have that for $j \ge 1$
  $$
  \Delta v_{j} - v_{j-1} = [\cdot]_{j-1} + \Delta u_{j}- [\cdot]_{j-1}- u_{j-1}
	=\Delta u_{j}- u_{j-1}
  \in \Pi_{j-2}
  $$
  and since $\left\{ v_0,\ldots,v_{j-2} \right\}$ is a basis of
  $\Pi_{j-2}$, the polynomial $\Delta v_{j} - v_{j-1}$ can be uniquely
  written as
  $$
   c_0 v_0+\cdots +c_{j-2} v_{j-2} =
  -\sum_{\ell=d-j+2}^{d} t_{d-j,\ell} \, v_{d-\ell}
  $$
  which defines the remaining entries of row $d-j$ of $\widetilde T_d$
  in a unique way such that $\widetilde T_d \vb=0$.
\end{pf}

\noindent
The last observation in the above proof can be formalized as follows.

\begin{corollary}\label{C:DualTaylor}
  For each $\vb \in \VV_d$ there exists a \emph{unique} generalized complete
  Taylor operator $\widetilde T_d$ such that $\widetilde T_d \vb = 0$. 
\end{corollary}

\begin{definition}
  The generalized Taylor operator of Corollary~\ref{C:DualTaylor},
  uniquely defined by
  \begin{equation}
    \label{eq:Annihil}
    \widetilde T (\vb) \, \vb = 0,
  \end{equation}
  is called the \emph{annihilator} of $\vb \in \VV_d$ and written as
  $\widetilde T (\vb)$. We 
	can skip 
	the subscript ``$d$'' because it
  is directly given by the dimension of $\vb$.
\end{definition}

\begin{definition}\label{D:ChainCompat}
  A \emph{chain} of length $d+1$ is a finite sequence $\Vb := [
  \vb_0,\dots,\vb_d ]$ of
  vectors
  $$
  \vb_j = \begin{bmatrix} 
      v_{j,j} \\ \vdots \\ v_{j,0}
    \end{bmatrix} =
    \begin{bmatrix}
      {[\cdot]_j} \\ \vdots \\ {[\cdot]_0}
    \end{bmatrix}
    + \ub_j \in \VV_j, \qquad j=0,\dots,d,
  $$
  that satisfies the \emph{compatibility condition}
  \begin{equation}
    \label{eq:ChainCompat}
    \wb_{j+1} :=
    \begin{bmatrix}
      w_{j+1,1} \\ \vdots \\ w_{j+1,j+1} \\
    \end{bmatrix}
    := \widetilde T (\vb_{j}) \begin{bmatrix} 
        v_{j+1,j+1} \\ \vdots \\ v_{j+1,1}
      \end{bmatrix} \in \RR^{j+1}, \qquad j=0,\dots,d-1.
    \end{equation}
\end{definition}

\begin{remark}
  Compatibility is a strong requirement on the interaction between
  $\vb_j$ and $\vb_{j+1}$. In general, $\widetilde T (\vb_{j}) \begin{bmatrix} 
    v_{j+1,j+1} \\ \vdots \\ v_{j+1,1}
  \end{bmatrix}$
  can only be expected to be a
  vector of polynomials in $\Pi_j,\dots,\Pi_0$, while compatibility
  requires all these polynomials to be constants.
\end{remark}

\noindent
Due to and by means of the compatibility condition, chains uniquely define a
generalized Taylor operator.

\begin{lemma}\label{L:ChainTaylwjj=1}
  If $\Vb$ is a chain of length $d+1$, then $w_{jj} = 1$, $j=1,\dots,d$.
\end{lemma}

\begin{pf}
  Since $v_{j+1,1} = [\cdot]_1 + c$ for some constant $c$ due to
  $\vb_j \in \VV_j$, it follows immediately from the definition
  (\ref{eq:ChainCompat}) that
  $$
  w_{j+1,j+1} = \Delta v_{j+1,1} = 1,
  $$
  as claimed.
\end{pf}

\begin{proposition}\label{L:ChainTaylProp}
  For $\Vb$ of length $d+1$ the following statements are equivalent:
  \begin{enumerate}
  \item\label{it:LChainTaylProp1} $\Vb$ is a chain of length $d+1$.
  \item\label{it:LChainTaylProp2} For $j=1,\dots,d$, we have
    \begin{equation}
      \label{eq:ChainTaylProp}
      \widetilde T ( \vb_j ) = \begin{bmatrix} 
          \widetilde T (\vb_{j-1})  & -\wb_j \\
          & \Delta
        \end{bmatrix} = \begin{bmatrix} 
          \Delta & -w_{1,1} & \dots & -w_{j,1} \\
                 & \Delta & \ddots & \vdots \\
                 & & \ddots & -w_{j,j} \\
                 & & & \Delta
        \end{bmatrix} .
    \end{equation}
  \item\label{it:LChainTaylProp3}
    \begin{equation}
      \label{eq:ChainTaylProp2}
      \widetilde T (\vb_d) \begin{bmatrix} 
          \vb_j \\ \BZero_{d-j}
        \end{bmatrix}  = 0, \qquad j=0,\dots,d.
    \end{equation}
  \end{enumerate}
\end{proposition}

\begin{pf}
  To show that \ref{it:LChainTaylProp1}) $\Rightarrow$
  \ref{it:LChainTaylProp2}), we 
  note that again (\ref{eq:ChainCompat}) yields that
  \begin{eqnarray*}
    0 & = & \widetilde T_j (\vb_{j}) \begin{bmatrix} 
      v_{j+1,j+1} \\ \vdots \\ v_{j+1,1}
    \end{bmatrix} 
    - \wb_{j+1}
    = \left[ \widetilde T (\vb_j) \,|\, -\wb_{j+1} \right] \begin{bmatrix} 
      v_{j+1,j+1} \\ \vdots \\ v_{j+1,1} \\ 1
    \end{bmatrix}  \\
      & = & \left[\widetilde  T (\vb_j) \,|\, -\wb_{j+1} \right] \, \vb_{j+1},
  \end{eqnarray*}
  Since $\widetilde T (\vb_{j+1})$ is unique, we deduce that
  \begin{equation}
    \label{eq:ChainTaylPropPf1}
    \widetilde T (\vb_{j+1}) = \begin{bmatrix} 
      \widetilde T (\vb_j) & -\wb_{j+1} \\
      & \Delta 
    \end{bmatrix} , \qquad
    j=0,\dots,d-1,
  \end{equation}
  which directly yields (\ref{eq:ChainTaylProp}).
	
  For \ref{it:LChainTaylProp2}) $\Rightarrow$
  \ref{it:LChainTaylProp3}) we simply notice that
  $$
  \widetilde T(\vb_d) \begin{bmatrix} 
      \vb_j \\ \BZero_{d-j}
    \end{bmatrix}  = \begin{bmatrix} 
    \widetilde  T(\vb_j) & * \\
      \BZero & * 
    \end{bmatrix}  \begin{bmatrix} 
      \vb_j \\ \BZero_{d-j}
    \end{bmatrix} = \begin{bmatrix} 
    \widetilde  T(\vb_j) \vb_j \\ \BZero
    \end{bmatrix} =\BZero,
  $$	
  while for \ref{it:LChainTaylProp3}) $\Rightarrow$
  \ref{it:LChainTaylProp1}) we first observe for $j < d$ that
  $$
  \BZero = \widetilde T(\vb_d) \begin{bmatrix} 
    \vb_j \\ \BZero_{d-j}
  \end{bmatrix} =
  \begin{bmatrix}
    \widetilde T(\vb_d)_{0:j,0:j} \vb_j \\ \BZero
  \end{bmatrix}
  $$
  and the uniqueness of the annihilators from Corollary
  \ref{C:DualTaylor} yields that $\widetilde T(\vb_d)_{0:j,0:j} =
  \widetilde T(\vb_j)$. This, in turn, implies together with
  \eqref{eq:ChainTaylProp2} that
  $$
  \BZero = \widetilde T(\vb_d) \begin{bmatrix} 
    \vb_{j+1} \\ \BZero
  \end{bmatrix}
  =
  \begin{bmatrix}
    \widetilde T (\vb_j) & -\wb_{j+1} & * \\
    & \Delta & * \\
    & & *
  \end{bmatrix}
  \begin{bmatrix}
    v_{j+1,j+1} \\ \vdots \\ v_{j+1} \\ 1 \\ \BZero
  \end{bmatrix}
  =
  \begin{bmatrix}
    \widetilde T (\vb_j) \,
    \begin{bmatrix}
      v_{j+1,j+1} \\ \vdots \\ v_{j+1,1} 
    \end{bmatrix}
    - \wb_{j+1} \\
    0 \\ \BZero 
  \end{bmatrix}
  $$
  which is the compatibility identity \eqref{eq:ChainCompat}, hence
  $\Vb$ is a chain.
\end{pf}

\begin{definition}\label{def:uniqtaylor}
  The unique generalized Taylor operator $\widetilde T (\vb_d)$ for a
  chain $\Vb$ will be written as $\widetilde T (\Vb)$.
\end{definition}

\begin{remark}
  The operator $\widetilde T( \Vb )$ of a chain $\Vb$ depends only on
  the last element $\vb_d$ of the chain.
\end{remark}

\begin{example}\label{Ex:ChainEx}
  Let $p_j = [\cdot]_j + q_j$, $q_j \in \Pi_{j-1}$, $j=0,\dots,d$, be
  given. Then
  $$
  \vb_j = \begin{bmatrix} 
      p_j, \, p_j',\,\dots,\, p_j^{(j)}
    \end{bmatrix}^T
  $$
  is a chain for the classical complete Taylor operator
  \begin{equation}
    \label{eq:ChainclassicalTayl}\widetilde T_{C,d} := \begin{bmatrix} 
      \Delta & -1 & -1/2! & -1/3!&\dots & -1/d! \\
      &\Delta & -1 & -1/2! & \dots & -1/(d-1)! \\
      && \Delta & -1 &  & \vdots \\
      &&  & \ddots & \ddots &  \vdots\\
      && & & \Delta & -1 \\
      && & & & \Delta
    \end{bmatrix}.
  \end{equation}
  Similarly,
  $$
  \vb_j = \begin{bmatrix} 
      p_j,\, \Delta p_j,\, \dots,\, \Delta^j p_j
    \end{bmatrix}^T
  $$
  is a chain for the operator
  \begin{equation}
    \label{eq:ChainExTDelta}
    \widetilde T_{\Delta,d} := \begin{bmatrix} 
      \Delta & -1 &&\BZero \\
             & \ddots & \ddots \\
             & & \ddots & -1 \\
             & & & \Delta
    \end{bmatrix} .    
  \end{equation}
  Another interesting generalized Taylor operator is
  \begin{equation}
    \label{eq:ChainExTSpline}
    \widetilde T_{S,d} := \begin{bmatrix} 
      \Delta & -1 & \dots & -1 \\
             & \ddots & \ddots & \vdots \\
             & & \ddots & -1 \\
             & & & \Delta
    \end{bmatrix},
  \end{equation}
  whose chains, connected to B--splines, we will consider in
  Example~\ref{Ex:SplineChain} later.
\end{example}

\noindent
As a shorthand for the property (\ref{eq:ChainTaylProp2}) of
Lemma~\ref{L:ChainTaylProp} we write $\widetilde T_d \Vb = 0$. Then we
have the following result.

\begin{lemma}\label{L:ChainExists}
  For any generalized complete Taylor operator $\widetilde T_d$ there
  exists a chain $\Vb$ of length $d+1$ such that $\widetilde T_d \Vb = 0$.
\end{lemma}

\begin{pf}
  The construction of the chain $\Vb$ is carried out inductively. To that
  end, we 
	recall 
	that if $p\in\Pi$ is of the form $\Delta
  p=[\cdot]_k$ for some $k\in\NN$, then $p=[\cdot]_{k+1}+c$ with some
  $c\in\RR$.

  Next, let $\vb_j \in \VV_j$, $j=0,\dots,d$, be any solution of
  $$
  0 = \widetilde T_d \begin{bmatrix}
      \vb_j \\ \BZero_{d-j}
    \end{bmatrix} = 
	\begin{bmatrix}     \widetilde T_j & * \\ 0 & *     \end{bmatrix}
	\begin{bmatrix}    \vb_j \\ \BZero_{d-j}     \end{bmatrix},
  $$
  or, equivalently, of $\widetilde T_j \vb_j = 0$. Such a solution can be
  found by setting $v_{j0} = 1$ and then solving, recursively for
  $k=1,\dots,j$, the equation given by row $j-k$ of the Taylor operator,
  \begin{equation}
    \label{eq:LChainExistsPf1}
    0 = \Delta v_{j,k} - v_{j,k-1} + \sum_{\ell=0}^{k-2} t_{j-k,j-\ell} \, v_{j,\ell}.
  \end{equation}
  Equivalently, this can be written with respect to the basis
  $\left\{ [\cdot]_0,\ldots,[\cdot]_{k-1} \right\}$ and using
  $v_{j,k-1}=[\cdot]_{k-1}+u_{j,k-1}$, $u_{j,k-1}\in\Pi_{k-2}$, as
  $$
  0=\Delta v_{j,k} - [\cdot]_{k-1}+ \sum_{\ell=0}^{k-2} s_{j-k,\ell} \,
  [\cdot]_\ell, \qquad s_{j-k,\ell} \in \RR,
  $$
  yielding
  $$
  v_{jk} = [\cdot]_k + \sum_{\ell=1}^{k-1} s_{j-k,\ell-1} [\cdot]_\ell +
  c_{k0}, \qquad k=0,\dots,j,
  $$
  where the constants $c_{k0} \in \RR$ can be chosen freely. This
  process yields polynomial vectors $\vb_j \in \VV_j$ such that
  $\widetilde T_j \vb_j = 0$, $j=0,\dots,d$. 
	
  Thus, it follows from the uniqueness of the annihilating Taylor
  operator from Corollary~\ref{C:DualTaylor} that $\widetilde T_j =
  T(\vb_j)$, and decomposing the identity
  $$
  0 = \widetilde T (\vb_{j+1}) \vb_{j+1}
  = \widetilde T_{j+1} \vb_{j+1} = \left[
    \begin{array}{cc}
      \widetilde T(\vb_j) & -\wb \\
      0 & \Delta
    \end{array}
  \right] \vb_{j+1}
  $$
  yields
  \begin{equation}
    \label{eq:LChainExistsPf2}
    \widetilde T (\vb_j) \left[
      \begin{array}{c}
        v_{j+1,j+1} \\ \vdots \\ v_{j+1,1}
      \end{array}
    \right] = \wb =: \wb_{j+1},    
  \end{equation}
  which is exactly the compatibility condition (\ref{eq:ChainCompat})
  needed for $\Vb$ to be a chain.
\end{pf}

\begin{corollary}\label{C:ChainConst}
  In the chain $\Vb$ from Lemma~\ref{L:ChainExists} the constant
  coefficients of the polynomials $v_{jk}$, $j=1,\dots,d$,
  $k=1,\dots,j$, can be chosen arbitrarily.
\end{corollary}

\begin{remark}
  The chain associated to a generalized Taylor operator is
  not at all unique, see also Example~\ref{Ex:ChainEx}.
\end{remark}

\noindent
The next result shows that any polynomial vector in $\VV_d$ can be
reached by a chain of length $d+1$.

\begin{proposition}
  For any $\vb \in \VV_d$ there exists a chain $\Vb$ of length $d+1$
  with $\vb_d = \vb$.
\end{proposition}

\begin{pf}
  Again we prove the claim by induction on $d$. The case $d=0$ is trivial as the
  only chain of length $0$ consists of $v=1$. For the induction step,
  we choose $\vb \in \VV_d$, $d > 0$ and the associated generalized
  Taylor operator $\widetilde T(\vb)$ as in Definition
  \ref{def:uniqtaylor}. 
  Then we know from   Lemma~\ref{L:ChainExists} that there exists a
  chain $\Vb = \left[ \vb_0,\dots,\vb_d \right]$ of length
  $d+1$ such that $\widetilde T (\vb) \Vb = 0$. Suppose that $\vb_d
  \neq \vb$ and, in particular, that $v_{d,1} (0) = v_1 (0) - 1$,
  which is possible according to Corollary~\ref{C:ChainConst}. With
  $$
  \vb = \left[
    \begin{array}{c}
      {[\cdot]_d} \\ \vdots \\ {[\cdot]_0}
    \end{array}
  \right] + \ub, \qquad \vb_d = \left[
    \begin{array}{c}
      {[\cdot]_d} \\ \vdots \\ {[\cdot]_0}
    \end{array}
  \right] + \ub_d, \qquad u_0 = u_{d,0} = 0,
  $$
  we find that
  $$
  0 = \widetilde T(\vb) \left( \vb - \vb_d \right) = \widetilde T(\vb) \left[
    \begin{array}{c}
      u_d - u_{d,d} \\ \vdots \\ u_1 - u_{d,1} \\ 0
    \end{array}
  \right] =: \widetilde T(\vb) \left[
    \begin{array}{c}
      \vb' \\ 0
    \end{array}
  \right]
  $$
  where $u_1 - u_{d,1} = v_1 (0) - v_{d,1} (0) = 1$. In addition,
  Lemma~\ref{L:TaylorKernels} yields that $\vb' \in \VV_{d-1}$ and
  therefore the decomposition
  $$
  \widetilde T(\vb) = \left[
    \begin{array}{cc}
      \widetilde T (\vb') & -\wb \\
      0 & \Delta
    \end{array}
  \right]
  $$
  and
  $$
  0 = \widetilde T(\vb) \vb = \left[
    \begin{array}{cc}
      \widetilde T (\vb') & -\wb \\
      0 & \Delta
    \end{array}
  \right] \left[
    \begin{array}{c}
      v_d \\ \vdots \\ v_1 \\ 1
    \end{array}
  \right] = \widetilde T (\vb') \left[
    \begin{array}{c}
      v_d \\ \vdots \\ v_1 \\ 1
    \end{array}
  \right] - \wb
  $$ compatibility between $\vb'$ and $\vb$. By the induction
  hypothesis, there exists a chain $\Vb'$ of length $d$ with $\vb_{d-1} =
  \vb'$ and since $\vb'$ is compatible with $\vb$, this chain can be
  extended to length $d+1$ with $\vb_d' = \vb$.
\end{pf}

\section{Chains and factorizations}
\label{sec:chainsfakt}

We now relate the existence of a spectral chain to factorizations of
the subdivision operators, thus extending the results first given in
\cite{MerrienSauer2012} for the classical Taylor operator.

\begin{definition}\label{D:SpecChain}
  A chain $\Vb$ of length $d+1$ is called \emph{spectral chain} for a
 vector subdivision scheme with mask $\Ab \in \ell^{(d+1)
    \times (d+1)} (\ZZ)$ if
  \begin{equation}
    \label{eq:SpecChain}
    S_\Ab \hat \vb_j = 2^{-j} \hat \vb_j, \qquad \hat \vb_j := \left[
      \begin{array}{c}
        \vb_j \\ \BZero_{d-j}
      \end{array}
    \right], \qquad j=0,\dots,d.
  \end{equation}
\end{definition}

\noindent
We will prove in Theorem \ref{T:SpecChainFact} 
that the existence of spectral chains is equivalent to the existence of
generalized Taylor factorizations. The main tool for this proof is the
following result.

\begin{proposition}\label{P:SCFact}
  If $\Cb \in 
  \ell^{(d+1) \times (d+1)} (\ZZ)$ 
  is a finitely supported
  mask for which there exists a chain $\Vb$ such that $S_\Cb \hat
  \vb_j = 0$, $j=0,\dots,d$, then there exists a finitely supported
  mask $\Bb \in 
  \ell^{(d+1) \times (d+1)} (\ZZ)$ 
  such that $S_\Cb
  = S_\Bb \, \widetilde T(\Vb)$.
\end{proposition}

\begin{pf}
  We follow the idea from \cite{MerrienSauer2012} and prove by
  induction on $k$ that the symbol $\Cb^* (z)$ satisfies
  \begin{equation}
    \label{eq:CBPartialDecomp}
    \Cb^* (z) =\Bb_k^* (z) \begin{bmatrix}
      \widetilde T(\vb_k)^* (z) & \BZero \\
      \BZero & \Ib
    \end{bmatrix}, \qquad k=0,\dots,d.
  \end{equation}
  with the columnwise written matrix
  \begin{equation}
    \label{eq:PSCFactPf1}
    \Bb_k^* (z) = \left[ \bb_0^* (z) \cdots \bb_k^* (z) \cb_{k+1}^* (z)
      \cdots \cb_d^* (z) \right].
  \end{equation}
  The construction makes repeated use of the well known factorization for a
  scalar subdivision scheme $S_a$:
  \begin{equation}
    \label{eq:factorscalar}
    \sum_{\alpha\in\ZZ}a(\alpha-2\beta)=0
    \qquad \Rightarrow \qquad a^*(z)=(z^{-2}-1)b^*(z),
  \end{equation}
  see, for example, \cite{CDM91} for proof.
	 
  For  case $k=0$, the annihilation of the vector $\hat \vb_0 =
  \eb_0=[1,0,\ldots,0]^T$ immediately gives the decomposition $\cb_0^*
  (z) = \left(z^{-2} - 1 \right) \bb_0^* (z)$ and therefore
  \begin{eqnarray*}
    \Cb^* (z)
    & = & \left[ \bb_0^* (z) \cb_1^* (z) \cdots \cb_d^* (z)
          \right] 
          \begin{bmatrix}
            z^{-2} - 1 \\
            & \Ib
          \end{bmatrix}\\
    & = & \left[ \bb_0^* (z) \cb_1^* (z) \cdots \cb_d^* (z)
          \right] \begin{bmatrix}
            \widetilde T(\vb_0)^* (z) \\
            & \Ib
          \end{bmatrix}.  
  \end{eqnarray*}
  Now suppose that (\ref{eq:CBPartialDecomp})
  holds for some $k \ge 0$. Then the fact that $\Vb$ is a chain
  yields, by means of the compatibility condition
  $$
  \wb_{k+1} = \widetilde T(\vb_k)
  \begin{bmatrix}v_{k+1,k+1}\\ \vdots \\v_{k+1,1}\end{bmatrix}
  $$
  that
  $$
  0 = S_\Cb \hat \vb_{k+1} =  S_{\Bb_k}
  \left[
    \begin{array}{cc|c}
      \widetilde T (\vb_k) & & \\
                           & 1 & \\
              \hline
                           & & \Ib 
    \end{array}
  \right] \begin{bmatrix}
      \vb_{k+1} \\ \BZero
    \end{bmatrix}
  = S_{\Bb_k} \begin{bmatrix}
      \wb_{k+1} \\ 1 \\ \BZero
   \end{bmatrix},
  $$
  or, applying \eqref{eq:factorscalar} to each row of the preceding equation, 
  $$
  \left[ \bb_0^* (z) \cdots \bb_k^* (z) \right]^T 
	\wb_{k+1}
	+ \cb_{k+1}^*
  (z) = \left( z^{-2} - 1 \right) \bb_{k+1}^* (z),
  $$
  which is
  $$
  \cb_{k+1}^* (z) = \left[ \bb_0^* (z) \cdots \bb_{k+1}^* (z) \right]^T 
	\begin{bmatrix}
      -\wb_{k+1} \\ z^{-2} - 1
    \end{bmatrix},
  $$
  or
  \begin{equation}
    \label{eq:PSCFactPf2}
    \Cb^* (z) = \left[ \bb_0^* (z) \cdots \bb_{k+1}^* (z) \cb_{k+2}^*
      (z) \cdots \cb_d^* (z) \right] \begin{bmatrix}
        \widetilde T(\vb_k)^* (z) & -\wb_{k+1} & \\
                                  & z^{-2} - 1 & \\
                                  & & \Ib
      \end{bmatrix}.    
  \end{equation}
  Since
  $$
  \widetilde T (\vb_{k+1})^* (z) =
  \begin{bmatrix}
      \widetilde T(\vb_k)^* (z) & -\wb_{k+1} \\
                                & z^{-2} - 1 \\
    \end{bmatrix},
  $$
  (\ref{eq:PSCFactPf2}) yields (\ref{eq:CBPartialDecomp}) with $k$
  replaced by $k+1$ and advances the induction hypothesis.
\end{pf}

\begin{remark}
  Proposition~\ref{P:SCFact} shows that, in the terminology of
  \cite{ContiCotroneiSauer17}, the generalized Taylor operator is a
  \emph{minimal annihilator} for the chain $\Vb$ since it annihilates
  the chain and factors any subdivision operator that does so, too.
\end{remark}

\noindent
Now we can show that the existence of a spectral chain results in the
existence of a factorization by means of generalized Taylor
operators. Since the Taylor operator corresponds to computing
differences, the scheme $S_\Bb$ from \eqref{eq:SpecChainFact} is often
called the \emph{difference scheme} of $S_\Ab$ with respect to the
generalized Taylor operator $\widetilde T (\Vb)$.

\begin{theorem}\label{T:SpecChainFact}
  If $S_\Ab$ possesses a spectral chain $\Vb$ of length $d+1$ then there
  exists a finite mask
  $\Bb \in \ell^{(d+1) \times (d+1)} (\ZZ)$ such that
  \begin{equation}
    \label{eq:SpecChainFact}
    \widetilde T(\Vb) \, S_\Ab = S_\Bb \, \widetilde T(\Vb).
  \end{equation}
\end{theorem}

\begin{pf}
  Since $S_\Cb := \widetilde T(\Vb) S_\Ab$ has the property that
  $$
  S_\Cb \hat \vb_k = \widetilde T(\Vb) S_\Ab \vb_k = 2^{-k} \widetilde
  T(\Vb) \vb_k = 0,
  $$
  an application of Proposition~\ref{P:SCFact} proves the claim.
\end{pf}

\begin{remark}
  For the validity of Theorem~\ref{T:SpecChainFact}, which is of a purely
  algebraic nature, the concrete eigenvalues of the spectral set are
  irrelevant. Their normalization will play a role, however, as soon
  as convergence is concerned. 
\end{remark}

\noindent
Next, we generalize a result from \cite{MerrienSauer2017} that serves
as a converse of Theorem~\ref{T:SpecChainFact}. The proof is a
modification of the former.

\begin{theorem}\label{T:FactorSpect}
  Suppose that for a finitely supported mask $\Ab \in \ell^{(d+1)\times (d+1)}$
  there exists a finitely supported $\Bb$ and a generalized
  incomplete Taylor operator $T_d$ such that $T_d S_\Ab = 2^{-d} S_\Bb T_d$ and 
  $S_\Bb \eb_d = \eb_d$. If a chain $\Vb$ for $T_d$ satisfies
  \begin{equation}
    \label{eq:FactorSpect}
    S_\Ab \hat \vb_j \in \mbox{\rm span} \left\{ \hat \vb_0,\dots,
      \hat \vb_j \right\}, \qquad j=0,\dots,d, 
  \end{equation}
  then there exists a spectral chain $\Vb'$ for $S_\Ab$.
\end{theorem}

\begin{pf}
  Relying on Lemma~\ref{L:ChainExists}, we choose a chain $\Vb$ such that
  $\widetilde T_d = \widetilde T (\Vb)$, which particularly yields
  that $T_d \vb_d = \eb_d$. Then
  $$
  T_d \vb_d = \eb_d = S_\Bb \eb_d = S_\Bb T_d \vb_d = 2^d T_d S_\Ab
  \vb_d
  $$
  implies that $T_d \left( 2^{-d} \vb_d - S_\Ab \vb_d \right) = 0$,
  hence
  $$
  S_\Ab \vb_d = 2^{-d} \vb_d + \tilde \vb, \qquad 0 = T_d \tilde \vb =
  \begin{bmatrix}
    \widetilde T_{d-1} & * \\ & 1
  \end{bmatrix} \tilde \vb,
  $$
  so that $\tilde v_0 = 0$ and therefore $\widehat T_{d-1} \tilde \vb_{0:d-1}
  = 0$.
  Since $\hat \vb_0,\dots,\hat \vb_{d-1}$ form a basis for the kernel of
  $\widetilde T_d$ with last component equal to zero, it follows that
  $\hat \vb \in \mbox{\rm span} \{ \vb_0,\dots,\vb_{d-1}\}$.
  Making use of the two--slantedness of
  $S_\Ab$, one can literally repeat the arguments of the proof of
  \cite[Theorem~2.11]{MerrienSauer2017} to conclude that
  $$
  S_\Ab \hat\vb_j - 2^{-j} \hat\vb_j \in \mbox{\rm span} \left\{
    \hat\vb_0,\dots,\hat\vb_{j-1} \right\},
  $$
  hence $S_\Ab \left[ \hat\vb_0,\dots,\hat\vb_d \right] = \left[
    \hat\vb_0,\dots,\hat\vb_d \right] \Ub$, where $\Ub \in \RR^{(d+1)
    \times (d+1)}$ is an upper triangular matrix with diagonal entries
  $1,\dots,2^{-d}$. Using the upper triangular $\Sb$ such that
  $\Sb^{-1} \Ub \Sb$ is diagonal, we can then define $\Vb'$ by
  $ \left[ \hat\vb_0',\dots,\hat\vb_d' \right] =
  \left[\hat\vb_0,\dots,\hat\vb_d \right] \Sb$,
  which is a chain since
  $$
  \widetilde T (\vb_d) \left( \sum_{k=0}^j c_k \, \hat\vb_k \right) = 0, \qquad
  j=0,\dots,d,
  $$
  due to Proposition~\ref{L:ChainTaylProp}.
\end{pf}

\section{Convergence}
\label{sec:Convergence}
From \cite{MerrienSauer2012,MerrienSauer2017}, we know that the
Hermite subdivision scheme  
$H_\Ab$ converges to a $C^d$ function according to
Definition~\ref{def:HermCdConv} if
\begin{enumerate}
\item there exists a scheme $S_\Bb$ such that $T_{C,d} S_\Ab = 2^{-d}
  S_\Bb T_{C,d}$ and $S_\Bb$ is convergent with limit function
  $\fb_\cb = \eb_d f_\cb$, where $\eb_d =[0,\ldots,0,1]^T$,
\item there exists a scheme $S_{\widetilde \Bb}$ such that $\widetilde
  T_{C,d} S_\Ab = 2^{-d} S_{\widetilde \Bb} \widetilde T_{C,d}$ and
  $S_{\widetilde \Bb}$ is \emph{contractive}. 
\end{enumerate}
Note that the normalization with the factor $2^{-d}$ now becomes
relevant since it affects the normalization and contractivity property
of $S_\Bb$ and $S_{\widetilde \Bb}$, respectively.


Before we give the results about the convergence replacing $T_{C,d}$
and $\widetilde T_{C,d}$ by $T$ and 
$\widetilde T$, respectively,
we will now consider conditions to guarantee that $\widetilde \Bb$ is
the result of such a factorization. To that end,
we recall  the factorization identity
\begin{equation}
  \label{eq:BBtildrelat}
  \begin{bmatrix}
      \Ib_d \\
   & \Delta
    \end{bmatrix}
   S_\Bb = S_{\widetilde \Bb} \begin{bmatrix}
      \Ib_d \\
   & \Delta
    \end{bmatrix}
\end{equation}
from vector subdivision \cite{MicchelliSauer98}. This relationship
does not depend on the form of the 
Taylor operator. In terms of symbols, (\ref{eq:BBtildrelat}) becomes
\begin{equation}
  \label{eq:BBtildrelatSymb}
  \begin{bmatrix}
      \Ib_d & \\
            & z^{-1} - 1
    \end{bmatrix}
 	\begin{bmatrix}
      \Bb_{11}^* (z) & \Bb_{12}^* (z) \\
      \Bb_{21}^* (z) & \Bb_{22}^* (z)
    \end{bmatrix} =
	\begin{bmatrix}
      \widetilde \Bb_{11}^* (z) & \widetilde \Bb_{12}^* (z) \\
      \widetilde \Bb_{21}^* (z) & \widetilde \Bb_{22}^* (z)
    \end{bmatrix}
	\begin{bmatrix}
      \Ib_d & \\
            & z^{-2} - 1
    \end{bmatrix}, 
\end{equation}
hence
\begin{eqnarray}
  \nonumber
  \Bb^* (z)
  & = & \begin{bmatrix} \Ib_d & \\ & z^{-1} - 1 \\ \end{bmatrix} ^{-1} 
  \begin{bmatrix}
    \widetilde \Bb_{11}^* (z) & \widetilde \Bb_{12}^* (z) \\
    \widetilde \Bb_{21}^* (z) & \widetilde \Bb_{22}^* (z) \\
  \end{bmatrix}
  \begin{bmatrix}\Ib_d & \\& z^{-2} - 1 \\
  \end{bmatrix} \label{eq:B*Btild*relat}\\
  & = & \begin{bmatrix}
    \widetilde \Bb_{11}^* (z)
    & (z^{-2} - 1) \, \widetilde \Bb_{12}^* (z) \\
    ( z^{-1} - 1 )^{-1} \widetilde \Bb_{21}^* (z)
    & (z+1) \, \widetilde \Bb_{22}^* (z) 
  \end{bmatrix}.
\end{eqnarray}

\begin{lemma}\label{L:BtildReqs}
  $S_\Bb$ converges to a continuous limit function of the form
  $\fb_\cb = f_\cb \, \eb_d$ if and only if
  $S_{\widetilde \Bb}$ is contractive,
  $\widetilde \Bb_{21} (1) = 0$ and
  $\widetilde \Bb_{22} (1) = 1$.
\end{lemma}

\begin{pf}
  First, to make $\Bb^*$ a Laurent polynomial, we must have $\widetilde
  \Bb_{21}^* (1) = 0$, otherwise $( z^{-1} - 1 )^{-1} \widetilde
  \Bb_{21}^* (z)$ has a pole at $1$. Second, the condition $S_\Bb
  \eb_d = \eb_d$ is equivalent to $\Bb^* (-1) \eb_d = 0$ and $\Bb^*
  (1) \eb_d = 2 \eb_d$. The first one of these requirements is
  automatically satisfied according to (\ref{eq:B*Btild*relat}), the
  second one becomes $2 \Bb_{22}^* (1) = 2$.
\end{pf}

Now we study the convergence of the Hermite scheme whenever we have one of the factorizations: 
$\widetilde TS_\Ab = 2^{-d}S_{\widetilde \Bb} \widetilde T$ or $T S_\Ab = 2^{-d}S_{\Bb}T$.
 To that end, we first recall the one dimensional case
of Lemma 3 in \cite{MerrienSauer2012}.

\begin{lemma}
 \label{L:conv1D}
  Given a sequence of refinements 
  $
  \hb_n=  \begin{bmatrix}
      h_n^{(0)} \\ h^{(1)}_n
    \end{bmatrix}
 \in\ell(\ZZ,\RR^{2})$ such that
  \begin{enumerate}
  \item 
    there exists a constant $c$ in $\RR$ such that 
    $\lim_{n\to +\infty} 
    h_n^{(0)} (0) =c$,
  \item
    there exists a  function $\xi \in C \left(\RR, \RR \right)$
    such that for any compact subset  $K$ of $\RR$ there exists a sequence
    $\mu_n$ 
    with limit $0$ and
    \begin{eqnarray}
      \label{eq:hyp2}
      \max_{\alpha \in 2^nK \cap \ZZ} \left | h_n^{(1)} (\alpha)- \xi
        \left( 2^{-n} \alpha \right)\right |_\infty & \le & \mu_n, \\
      \label{eq:hyp3}
      \max_{\alpha \in 2^nK \cap \ZZ} \left | 2^n \Delta h^{(0)}_n (\alpha) -
        h_n^{(1)} (\alpha)\right |_\infty & \le & \mu_n.
    \end{eqnarray}
  \end{enumerate}
  Then there exists for any compact $K$ 
  a sequence $\theta_n$ with limit $0$ such that 
  the function
  \begin{equation}
    \label{eq:phiintDef}
    \varphi(x) = c + \int_0^1 x \, \xi \left( tx \right) dt, \qquad x \in
    \RR,    
  \end{equation}
  satisfies
  \begin{equation}
    \label{eq:majorfin}
    \max_{\alpha\in2^nK\cap\ZZ} \left\| h_n^{(0)} (\alpha) - \varphi \left(
        2^{-n} \alpha \right) \right\| \le \theta_n, \qquad n \in \NN.
  \end{equation}
\end{lemma}

\begin{theorem}
  \label{T:convdD}
  Let $\Ab,\Bb \in \ell^{d+1}(\ZZ)$ be two masks related by the 
  the factorization $T_d S_\Ab=2^{-d}S_\Bb T_d$ for some generalized
  incomplete Taylor operator $T_d$.
	
  Suppose that for any initial data 
  $\fb_0 \in \ell^{d+1} ( \ZZ)$ and associated refinement
  sequence $\fb_n$ of the Hermite scheme $H_\Ab$,
  \begin{enumerate}
  \item
    the sequence 
    $\fb_n (0)$ converges to a limit $\yb \in\RR^{d+1}$,
  \item
    the subdivision scheme $S_\Bb$ is 
    $C^{p-d}$--convergent for some $p\ge d$, and that for any initial
    data $\gb_0 = T_d \fb_0$, the limit function $\Psi = \Psi_\gb \in
    C^{p-d} \left(\RR,\RR^{d+1} \right)$ satisfies 
    \begin{equation}
      \label{eq:PsiProp}
      \Psi = \begin{bmatrix}
          \BZero \\
          \psi 
        \end{bmatrix}, \qquad
        \psi\in C^{p-d}\left(\RR,\RR \right).
    \end{equation}
  \end{enumerate}
  Then $H_\Ab$ is $C^p$--convergent.
\end{theorem}

\begin{pf} The proof is adapted from the proofs in
  \cite{DubucMerrien09,MerrienSauer2011}. Given $\fb_0 \in
  \ell^{d+1}(\ZZ)$, let $\gb_0 = T_d \fb_0$. We define the following
  two sequences: $\fb_{n+1}=\Db^{-n-1}S_\Ab
  (\Db\fb_n)$ and $\gb_{n+1}=S_\Bb\gb_n$, $n \in \NN$.  
  Since $T_d S_\Ab=2^{-d}S_\Bb T_d$, we can directly deduce that
  $\gb_n =2^{nd} T_d \Db^n \fb_n$. 
	
  With the convergence of $\fb_n (0)$, let $y_i:=\lim_{n\to
    +\infty}f_n^{(i)} (0)$, $i=0,\ldots,d$.
  Then we define $\Phi$ recursively beginning with $\phi_d = \psi$
  and setting 
  \begin{equation}
    \label{def_phi_a}
    \phi_i(x) = y_i+ \int_0^1 x \, \phi_{i+1}(tx) \, dt \qquad i =
    d-1,\dots,0.
  \end{equation}
  Then $\Phi=[\phi_i]_{i=0,\ldots d}$ is continuous with
  $\phi_i^{(d-i)}=\psi$.

  Fixing a compact  $K \subset \RR$,
  we will prove by a backward finite recursion that for $k=d,d-1,\dots,0$, 
	 there exists a sequence $\varepsilon_n$ with limit $0$ such that
  \begin{equation}
    \label{propertypk}
    \left| f_n^{(k)} (\gamma) - \phi_k\left( 2^{-n} \gamma \right)
    \right| \le \varepsilon_n, \qquad \gamma\in \ZZ\cap 2^nK.
  \end{equation}
  The case $k=d$ is an immediate consequence of the convergence of the last row of $\gb_n$ and $g_n^{(d)}=f_n^{(d)}$, 
	which yields for any $\gamma\in\ZZ\cap2^nK$ that
  \begin{equation}
    \label{derivd}
    \left| f_n^{(d)} (\gamma) - \psi ( 2^{-n} \gamma ) \right|\le \varepsilon_n,
  \end{equation}
  while, for $k < d$,
  the convergence of 
  the appropriate component of $\gb_n$ 
  to zero implies that 
  \begin{equation}
    \label{derivk}
    2^{n(d-k)} \left|\Delta f_n^{(k)}( \gamma)-\frac 1{2^n}f_n^{(k+1)}( \gamma)
    +  \sum_{\ell=2}^{d-k} t_{k,k+\ell}\frac 1{2^{n\ell}} f_n^{(k + \ell)} ( \gamma)\right|\le \varepsilon_n,
  \end{equation}
  for a sequence $\varepsilon_n$ that tends to zero for $n \to \infty$.
  
  To prove \eqref{propertypk} for $k=d-1$, we define
  the sequences $\hb_n =[f_n^{(d-1)}, f_n^{(d)}]^T$. Then \eqref{derivk} becomes
	$\left|2^n\Delta f_n^{(d-1)}( \cdot)-f_n^{(d)}( \cdot)\right|\le \varepsilon_n$. 
	Because of \eqref{derivd}, we
  can apply Lemma \ref{L:conv1D} and obtain that
  $$
  \left| f_n^{(d-1)}( \gamma ) - \phi_{d-1} \left( 2^{-n} \gamma
    \right) \right| \le \theta_n, \qquad \gamma\in2^nK\cap\ZZ^,
  $$
  which is \eqref{propertypk} for $k=d-1$.
  
  To prove the recursive step $k+1 \to k$, $0 \le k < d-2$, we 
  get from \eqref{derivk} that, for $\gamma\in\ZZ\cap 2^nK$,
  \begin{equation}
    \label{eq:epsnkBound}
    \left| 2^n \Delta f_n^{(k)} (\gamma) - f_n^{(k+1)} (\gamma)
    \right| \le \frac{\varepsilon_n}{2^{n(d-k)}} + \sum_{\ell=2}^{d-k}
    \frac{|t_{k,k+\ell}|}{2^{n\ell}} \left| f_n^{(k+\ell)}(\gamma)\right|
  \end{equation}
  Since \eqref{propertypk} holds for $j > k$, it follows that
  $$
  \lim_{n \to \infty} \left| f_n^{(j)} (\gamma)-
    \phi_j \left( 2^{-n} \gamma \right) \right| = 0
  $$  
  uniformly for $\gamma\in\ZZ\cap 2^nK$ and since $\phi_j$ is bounded on $K$, so is the sequence $\left|
    f_n^{(j)} (\gamma) \right|$ on $\ZZ\cap 2^nK$. 
  Thus the right hand side of \eqref{eq:epsnkBound} 
  converges to zero so
  that it immediately implies \eqref{propertypk} using again Lemma \ref{L:conv1D}.
\end{pf}

\noindent
As a consequence of Theorem \ref{T:convdD} and Lemma \ref{L:BtildReqs}
we also have the following result.

\begin{corollary}
 \label{T:convdDcomplete}
  Let $\Ab,\tilde \Bb \in \ell^{d+1}(\ZZ)$ be two masks related by the 
  the factorization $\tilde T_d S_\Ab=2^{-d}S_{\tilde \Bb} \tilde T_d$
  for some generalized 
  complete Taylor operator $\tilde T_d$. For  any initial data 
  $\fb_0 \in \ell^{d+1} ( \ZZ)$ and associated refinement
  sequence $\fb_n$ of the Hermite scheme $H_\Ab$, we suppose that
    the sequence 
    $\fb_n (0)$ converges to a limit $\yb \in\RR^{d+1}$.
If  $S_{\widetilde \Bb}$ is contractive,  $\widetilde \Bb_{21} (1) = 0$ and
  $\widetilde \Bb_{22} (1) = 1$,  then $H_\Ab$ is $C^d$--convergent.
\end{corollary}

\begin{remark}
  The condition that $\fb_n (0)$ converges can be discarded by using
  the techniques from \cite{ccs17:_hermit}. The
  factorization arguments used there can easily be seen to carry over
  to the situation of arbitrary generalized Taylor
  operators. Nevertheless, we prefer the proof given here due to its
  analytic flavor which nicely corresponds to the graphs shown
  later. There the function $\psi$ equals the last derivative of the
  limit function in accordance with the proof above.
\end{remark}

\section{Unfactoring constructions}
\label{sec:unfactor}
In this section we consider the construction of convergent Hermite
subdivision schemes that factorize with respect to a given generalized
Taylor operator, thus showing that there exist whole classes of
convergent Hermite subdivision schemes that do \emph{not} satisfy the
spectral condition. In particular, the spectral condition is not
necessary for $C^d$--convergence.

These constructions will be based on determining a contractive
difference scheme $\widetilde \Bb$. The difficulty, as in all vector
subdivision schemes, lies in the fact that, in contrast to the scalar
case, not every vector subdivision scheme is the 
difference scheme of a finitely supported vector or Hermite
subdivision schemes, but that more intricate algebraic conditions have
to be taken into account.

\subsection{Conditions on the difference schemes}
We begin with an inversion of the Taylor operator.

\begin{lemma}\label{L:TaylorInverse}
  For any generalized complete Taylor operator $\widetilde T_d$, 
  there exists an upper triangular matrix $\Pb^* (z)$ of Laurent
  polynomials such that
  \begin{equation}
    \label{eq:TaylorInverse}
    \left( \widetilde 
      T_d^*
      (z) \right)^{-1} = \frac{1}{z^{-1}-1}
    \Db_d^* (z) \Pb^* (z) \left( \Db_d^* (z) \right)^{-1},
  \end{equation}
  where
  $$
  \Db_d^* (z) = \begin{bmatrix}
    1 \\
    & z^{-1} - 1 \\
    & & \ddots \\
    & & & ( z^{-1} - 1 )^d
  \end{bmatrix}.
  $$
  Moreover $p_{jj}^* (z) = 1$, $j=0,\dots,d$, and
  \begin{equation}
    \label{eq:P*1}
    \Pb^* (1) = 
    \begin{bmatrix}
      1 & \dots & 1 \\
      & \ddots & \vdots \\
      & & 1
    \end{bmatrix}.    
  \end{equation}
\end{lemma}

\begin{pf}
  Since
  \begin{eqnarray*}
    \widetilde 
     T_d^* (z)
    & = & \begin{bmatrix}
      z^{-1} - 1 & * & \dots & * \\
      & z^{-1} - 1 & \ddots & \vdots \\
      & & \ddots & * \\
      & & & z^{-1} - 1 \\
    \end{bmatrix}
             =  ( z^{-1} - 1 ) \left( \Ib -
             \frac{\Nb}{z^{-1} - 1} \right)
  \end{eqnarray*}
  with the strictly upper triangular nilpotent matrix
  $$
  \Nb = \begin{bmatrix}
    0 & 1 & * & \dots & * \\
    & 0 & 1 & \ddots & \vdots \\
    & & 0 & \ddots & * \\
    & & & \ddots & 1 \\
    & & & & 0
  \end{bmatrix}
  \in \RR^{(d+1) \times (d+1)}, \qquad \Nb^{d+1} = 0,
  $$
  it follows that
  \begin{eqnarray*}
    \lefteqn{\left( \widetilde 
          T_d^* (z) \right)^{-1}
    = \frac{1}{z^{-1} - 1}
    \left( \Ib + \sum_{j=1}^d \left( \frac{\Nb}{z^{-1} - 1} \right)^j
    \right)} \\
 & = & \begin{bmatrix}
   \frac{p_{00}^* (z)}{z^{-1} - 1} & \frac{p_{01}^* (z)}{(z^{-1}- 1)^2 } & \dots 
   & \frac{p_{0d}^* (z)}{(z^{-1} - 1)^{d+1}}
   \\ & \frac{p_{11}^* (z)}{z^{-1} - 1} & \ddots & \vdots \\ 
   & & \ddots & \frac{p_{d-1,d}^* (z)}{(z^{-1} -   1)^2 } \\  
   & & & \frac{p_{dd}^* (z)}{z^{-1} - 1}\end{bmatrix}
         = \frac{1}{z^{-1}-1} \Db_d^* (z) \Pb^* (z) \left( \Db_d^* (z)
         \right)^{-1}.
  \end{eqnarray*}  
  The property of the diagonal elements $p_{jj}$ is immediate from the form of
  $\Nb$, in particular $ \sum_{j=1}^d \left( \frac{\Nb}{z^{-1} - 1} \right)^j$ has a null diagonal.
	
  For the computation on the off-diagonal elements, we notice that due to
  $$
  \Nb^j = 
  \begin{bmatrix}
    0 & \dots & 0 & 1 & * & \dots & * \\
    & \ddots & & \ddots & \ddots & \ddots & \vdots \\
    & & 0 & \dots & 0 & 1 & * \\
    & & & \ddots & & \ddots & 1 \\
    & & & & 0 & \dots & 0 \\
    & & & & & \ddots & \vdots \\
    & & & & & & 0 \\
  \end{bmatrix},
  $$
  it follows that
  $$
  \frac{p_{jk}^* (z)}{(z^{-1} - 1)^{k-j+1}} = \frac{1}{(z^{-1} - 1)^{k-j+1}} +
  \frac{q_{jk} (z)}{(z^{-1} - 1)^{k-j}}
  = \frac{(z^{-1} - 1) q_{jk} (z) + 1}{(z^{-1} - 1)^{k-j+1}},
  $$
  which gives (\ref{eq:P*1}).
\end{pf}

\begin{example}
  For the generalized complete Taylor operator $\widetilde T_\Delta$
  from \eqref{eq:ChainExTDelta}, we get the constant polynomial
  $$
  \Pb^* (z) = \Pb^* (1) = \left[
    \begin{array}{ccc}
      1 & \dots & 1 \\
        & \ddots & \vdots \\
        & & 1
    \end{array}
  \right].
  $$
\end{example}

\noindent
Next, we compute $\left( \widetilde 
    T_d^* (z) \right)^{-1}
\widetilde \Bb^* (z)$, by first noting that
$$
\frac{1}{z^{-1} - 1} \left( \Db_d^* (z) \right)^{-1} \widetilde \Bb^*
(z) = \left[
  \begin{array}{ccc}
    \frac{\widetilde b_{00}^* (z)}{z^{-1} - 1} & \dots
    &\frac{\widetilde b_{0d}^*(z)}{z^{-1} -1}  \\ 
    \vdots & \ddots & \vdots \\
    \frac{\widetilde b_{d0}^* (z)}{(z^{-1} - 1)^{d+1}} & \dots
    &\frac{\widetilde b_{dd}^*(z)}{(z^{-1} -1)^{d+1}} 
  \end{array}
\right].
$$
Therefore the entries $c_{jk}^* (z)$ of
$$
\Cb^* (z) := \left( \widetilde T_d^* (z) \right)^{-1} \widetilde \Bb^*
(z)
= (z^{-1} - 1)^{-1} \Db_d^* (z) \, \Pb^* (z) \left( \Db_d^* (z)
\right)^{-1} \widetilde \Bb^* (z)
$$
satisfy, for $j,k=0,\dots,d$,
\begin{eqnarray*}
  c_{jk}^* (z)
  & = & (z-1)^j \sum_{\ell=j}^d p_{j{\ell}}^* (z) \frac{\widetilde b_{{\ell}k}^*
        (z)}{(z^{-1} - 1)^{\ell+1}}
        = \sum_{\ell=j}^d p_{j{\ell}}^* (z) \frac{\widetilde b_{{\ell}k}^*
        (z)}{(z^{-1} - 1)^{\ell-j+1}}.
\end{eqnarray*}
Then, the components $a_{jk}^* (z)$ of the final result
$$
\Ab^* (z) = \left( \left( \widetilde 
		T_d\right)^* (z) \right)^{-1}
\widetilde \Bb^* (z) \left( \widetilde 
		T_d\right)^* (z^2) =
\Cb^* (z) \, \left( \widetilde 	T_d\right)^* (z^2)
$$
satisfy, since $\left( \left( \widetilde 	T_d\right)^* (z^2)
\right)_{rk} = 0$ for $r > k$,
\begin{eqnarray*}
  a_{jk}^* (z)
  & = & \sum_{r=0}^d c_{jr}^* (z) \, \left( \left( \widetilde 
		T_d \right)^* (z^2) \right)_{rk} =
        \sum_{r=0}^k c_{jr}^* (z) \, \left( \left( \widetilde 
						T_d\right)^* (z^2) \right)_{rk} \\
  & = & (z^{-2} - 1) c_{jk}^* (z) - \sum_{r=0}^{k-1}  c_{jr}^* (z) \,
        w_{k,r+1} \\
  & = & (z^{-1} + 1 ) \sum_{\ell=j}^d p_{j{\ell}}^* (z) \frac{\widetilde b_{{\ell}k}^*
        (z)}{(z^{-1} - 1)^{\ell-j}} - \sum_{r=0}^{k-1} w_{k,r+1}
        \sum_{\ell=j}^d p_{j{\ell}}^* (z) \frac{\widetilde
        b_{{\ell}r}^* (z)}{(z^{-1} - 1)^{\ell-j+1}},
\end{eqnarray*}
hence,
\begin{equation}
  \label{eq:ajkFormula}
  a_{jk}^* (z) = \sum_{\ell=j}^d \frac{p_{j{\ell}}^* (z)}{(z^{-1} -
    1)^{\ell-j}} \left( (z^{-1} +  1 ) \widetilde b_{{\ell}k}^* (z) -
    \sum_{r=0}^{k-1} w_{k,r+1} \frac{\widetilde b_{{\ell}r}^* (z)}{z^{-1} - 1}
  \right), \qquad j,k=0,\dots,d.
\end{equation}

\begin{lemma}\label{L:blinecond}
  If for any $j,k=0,\ldots,d$, there exists a Laurent polynomial
  $h_{jk}^*(z)$ such that 
  \begin{equation}
    \label{eq:blinecond}
    (z^{-1} +  1 ) \widetilde b_{jk}^*
    (z) - \sum_{r=0}^{k-1} w_{k,r+1} \frac{\widetilde b_{jr}^*
        (z)}{z^{-1} - 1} 
								=
				(z^{-1} - 1 )^j h_{jk}^* (z), 
    \end{equation}
    then $\Ab \in \ell^{(d+1) \times (d+1)} (\ZZ)$.
\end{lemma}

\begin{pf}
  Since $p_{j\ell}^* (1) = 1$, all the terms of the outer sum in
  (\ref{eq:ajkFormula}) are
  polynomials if and only if
  $$
  (z^{-1} +  1 ) \widetilde b_{{\ell}k}^*
  (z) - \sum_{r=0}^{k-1} w_{k,r+1} \frac{\widetilde b_{{\ell}r}^*
    (z)}{z^{-1} - 1}, \qquad \ell=j,\dots,d,
  $$
  has an $(\ell-j)$--fold zero at $1$ for all $j \le \ell$, in particular for
  $j = 0$, which yields (\ref{eq:blinecond})
  after replacing $\ell$ by $j$.
 \end{pf}

\noindent
The simplest way to solve (\ref{eq:blinecond}) is to set
\begin{equation}
  \label{eq:bjkFactorSimple}
  \widetilde b_{jk}^* (z) = 
  (z^{-1}-1)^j h_{jk}^* (z), \qquad j=0,\dots,d-1, \quad k=0,\dots,d,  
\end{equation}
which we can even choose in a upper triangular way by setting
$h_{jk}^* = 0$ for $k > j$. Note that this choice is even independent
of the generalized Taylor operator.

For the final row, however, we cannot use this approach since it would yield
$\widetilde b_{dd}^* (1) = 0$, thus contradicting the requirement
from Lemma~\ref{L:BtildReqs}. To overcome this problem, we set
\begin{equation}
  \label{eq:hFactdef}
  \widetilde b_{dj}^* (z) = (z^{-1} - 1) \, g_{dj}^* (z) =: (z^{-1} -
  1)^{d-j} \, h_{dj}^* (z^{-1}), \qquad j=0,\dots,d.
\end{equation}
We want to construct $h_{dj}^*$ in such a way that for $j=0,\dots,d$
the polynomials
\begin{eqnarray*}
  \lefteqn{(z^{-1} +  1 ) \widetilde b_{dj}^*
  (z) - \sum_{k=0}^{j-1} w_{j,k+1} \frac{\widetilde b_{dk}^*
  (z)}{z^{-1} - 1}}\\
  & = & (z^{-1} +  1 ) ( z^{-1} - 1 )^{d-j} h_{dj}^* (z^{-1}) - \sum_{k=0}^{j-1}
        w_{j,k+1} (z^{-1} - 1 )^{d-k-1} h_{dk}^* (z^{-1}) \\
  & = & (z^{-1} - 1 )^{d-j} \left( (z^{-1}+1) h_{dj}^* (z^{-1}) -
        \sum_{k=0}^{j-1}  w_{j,k+1} \, (z^{-1}-1)^{(j-1) - k} h_{dk}^*
        (z^{-1}) \right) \\
  & = & (z^{-1} - 1 )^{d-j} \left( (z^{-1}+1) h_{dj}^* (z^{-1}) -
        \sum_{k=0}^{j-1}  w_{j,j-k} \, (z^{-1}-1)^k h_{d,j-1-k}^*
        (z^{-1}) \right)
\end{eqnarray*}
have a zero of order $d$ at $1$. Since $w_{jj} = 1$, this is
equivalent, after replacing $z$ by $z^{-1}$,
to a zero of order $j$ at $1$ of the Laurent polynomials
\begin{equation}
  \label{eq:hdjrecurrence}
  q_j (z) := (z+1) h_{dj}^* (z) - h_{d,j-1}^* (z) - \sum_{k=1}^{j-1} w_{j,j-k} \,
  (z-1)^k h_{d,j-1-k}^* (z).
\end{equation}
This implies that
$$
0 = q_j (1) = 2 h_{dj}^* (1) - h_{d,j-1}^* (1), \qquad j=1,\dots,d,
$$
which yields, together with the requirement that $\widetilde b_{dd}^*
(1) = 1$, that
\begin{equation}
  \label{eq:hcondat1}
  h_{dj}^* (1) = 2^{d-j}, \qquad j=0,\dots,d.
\end{equation}
The $r$th derivative, $r=1,\dots,j$, of $q_j$ is
\begin{eqnarray*}
  q_j^{(r)} (z)
  & = & \sum_{s=0}^r {r \choose s} \frac{d^s}{dz^s} (z+1) \, \left(
        h_{dj}^* \right)^{(r-s)} (z) - \left( h_{d,j-1}^* \right)^{(r)}
        (z) \\
  && \quad - \sum_{k=1}^{j-1} w_{j,j-k} \, \sum_{s=0}^r {r \choose s}
     \left( \frac{d^s}{dz^s} (z-1)^k \right) \, \left( h_{d,j-1-k}^*
     \right)^{(r-s)} (z) \\
  & = & (z+1) \left( h_{dj}^* \right)^{(r)} (z) + r \, \left( h_{dj}^*
        \right)^{(r-1)} (z) - \left( h_{d,j-1}^* \right)^{(r)}
        (z) \\
  && \quad - \sum_{k=1}^{j-1} w_{j,j-k} \, \sum_{s=0}^{\min (k,r)} {r \choose s}
     \frac{k!}{(k-s)!} (z-1)^{k-s}  \, \left( h_{d,j-1-k}^*
     \right)^{(r-s)} (z).
\end{eqnarray*}
Therefore, we can express the additional requirements as
\begin{eqnarray}
  \nonumber
  0
  & = & q_j^{(r)} (1) \\
  \nonumber
  & = & 2 \left( h_{dj}^* \right)^{(r)} (1) + r \, \left( h_{dj}^*
        \right)^{(r-1)} (1) - \left( h_{d,j-1}^* \right)^{(r)} (1) \\
  \label{eq:qDerivCondr}
  & & \quad - \sum_{k=1}^{r} w_{j,j-k} \, \frac{r!}{(r-k)!}
      \left( h_{d,j-1-k}^*
      \right)^{(r-k)} (1), \qquad r =1,\dots,j-1,
\end{eqnarray}
and, with $r = j$,
\begin{eqnarray}
  \nonumber
  0 & = & 2 \left( h_{dj}^* \right)^{(j)} (1) + r \, \left( h_{dj}^*
          \right)^{(j-1)} (1) - \left( h_{d,j-1}^* \right)^{(j)} (1)
  \\
  \label{eq:qDerivCondj}
  & & \quad - \sum_{k=1}^{j-1} w_{j,j-k} \, \frac{j!}{(j-k)!}
      \left( h_{d,j-1-k}^* \right)^{(j-k)} (1).
\end{eqnarray}
Together, \eqref{eq:qDerivCondr} and \eqref{eq:qDerivCondj} can be
used to build the polynomials $h_{dj}^*$ recursively.

This construction allows us to easily create factorizable schemes via
\eqref{eq:qDerivCondr} and \eqref{eq:qDerivCondj}, but it is more difficult
to choose $h_{d0}^* (z)$ in such a way that the final $h_{dd}^* (z)$
is the symbol of a contractive scheme. To achieve this, we perform the
recurrence in the opposite direction, which is still easy for $\widetilde
T_\Delta$.

\begin{figure}
  \centering
  \includegraphics[width=.75\hsize]{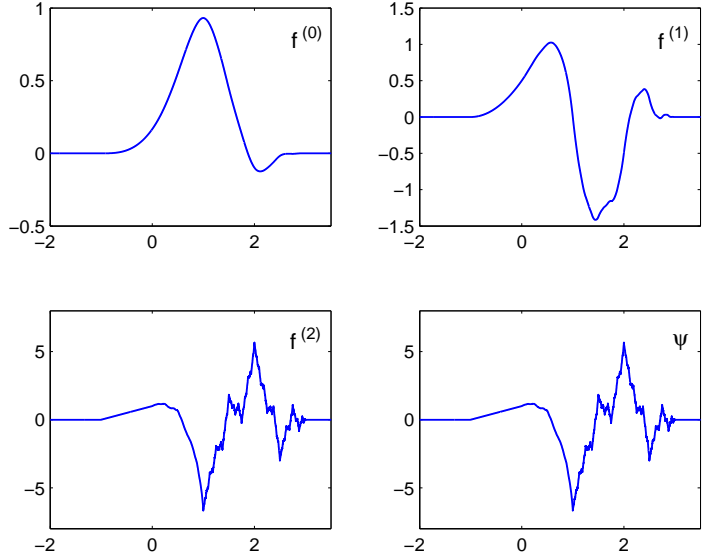}
  \caption{Limit functions for Example~\ref{Ex:DeltaEx}, showing the
    three entries of the limit function of the Hermite subdivision
    scheme and the limit function of the associated convergent
    difference scheme.}
  \label{fig:Exdelta}
\end{figure}

\begin{example}\label{Ex:DeltaEx}
  For the generalized Taylor operator $\widetilde T_\Delta$ we get
  the simplified conditions
  \begin{equation}
    \label{eq:DelthCond}
    0 = 2 \left( h_{d,j}^* \right)^{(r)} (1) + r \, \left( h_{d,j}^*
    \right)^{(r-1)} (1) - \left( h_{d,j-1}^* \right)^{(r)} (1), \qquad
    r=1,\dots,j,
  \end{equation}
  or
  \begin{equation}
    \label{eq:DelthCond2}
    \left( h_{d,j}^* \right)^{(r)} (1) = \frac12 \left( \left(
        h_{d,j-1}^* \right)^{(r)} (1) - r \, \left( h_{d,j}^*
      \right)^{(r-1)} (1) \right), \qquad r=1,\dots,j.
  \end{equation}
  To come up with convergent schemes of arbitrary size that factor
  with $\widetilde T_\Delta$, we now solve (\ref{eq:DelthCond})  for
  $h_{d,j-1}^*$, replace $j-1$ by $j$ and thus get
  $$
  \left( h_{dj}^* \right)^{(r)} (1) = 2 \left( h_{d,j+1}^*
  \right)^{(r)} (1) + r \, \left( h_{d,j+1}^* \right)^{(r-1)} (1), \qquad
  r = 1,\dots,j+1,
  $$
  which leads to the the explicit formula
  \begin{equation}
    \label{eq:exj->j-1rule}
    h_{dj}^* (z) = 2^{d-j} + \sum_{r=1}^{n+d-j} \frac{2 \left( h_{d,j+1}^*
      \right)^{(r)} (1) + r \left( h_{d,j+1}^* \right)^{(r-1)}
      (1)}{r!} \, (z-1)^r, \qquad j=d-1,\dots,0,
  \end{equation}
  initialized with a polynomial $h_{dd}^*$ of degree $n$. Starting
  with the simplest choice $h_{dd}^* (z) = \frac12 ( z+1)$, we thus
  get
  \begin{eqnarray*}
    h_{d,d-1}^* (z) & = & 2 + 2 ( z-1 ) + \frac12 (z-1)^2 = \frac12 (z
                          + 1)^2 \\
    h_{d,d-2}^* (z) & = & 4 + 6 (z-1) + 3 (z-1)^2 + \frac12 (z-1)^3 =
                          \frac12 (z + 1)^3.
  \end{eqnarray*}
  If we now set $f_n
  (z) := \frac12 (z+1)^n$, then $f_n^{(r)} (1) = \frac{n!}{(n-r)!}
  2^{n-r-1}$ and the fact that
  \begin{eqnarray*}
    \lefteqn{
    2 f_{n-1}^{(r)} (1) + r f_{n-1}^{(r-1)} (1) =
    \frac{(n-1)!}{(n-1-r)!} 2^{n-1-r} + r \, \frac{(n-1)!}{(n-r)!}
    2^{n-1-r}} \\
    & = & \frac{(n-1)!}{(n-1-r)!} 2^{n-1-r} \left( 1 + \frac{r}{n-r}
          \right)
          = \frac{(n-1)!}{(n-1-r)!} 2^{n-1-r} \frac{n}{n-r}
          = \frac{n!}{(n-r)!} 2^{n-r-1} \\
    & = & f_n^{(r)} (1)
  \end{eqnarray*}
  shows that indeed
  \begin{equation}
    \label{eq:TDelthGeneric}
    h_{d,j}^* (z) = \frac12 (z+1)^{d-j+1}, \qquad j=0,\dots,d,
  \end{equation}
  satisfy the recurrence \eqref{eq:exj->j-1rule} and therefore
  $$
  \widetilde b_{dj}^* (z) = (z^{-1}-1)^{d-j} \, h_{dj}^* (z^{-1}) =
  \frac12 \left( z^{-1} + 1 \right) \, \left( z^{-2} - 1 \right)^{d-j}
  $$
  is a proper choice. For $d=2$, for example, we can set
  $$
  \widetilde \Bb^* (z) =
  \begin{bmatrix}
    -{{z-1}\over{2 \, z}}&0&0\\
    {{\left(z-1\right)^2}\over{z^2}}&
    {{\left(z-1\right)^2}\over{4 \, z^2}}&0\\
    {{\left(z-1\right)^2\,\left(1+
          z\right)^3}\over{2\,z^5}}&-{{\left(z-1\right)\,\left(1+z\right)^2
      }\over{2\,z^3}}&{{1+z}\over{2\,z}}\\
  \end{bmatrix}
  $$
  and get the corresponding
  $$
  \Ab^* (z) =1/4
  \begin{bmatrix}
    -{{\left(1+z\right)\,\left(-1-3\,z-6\,z^2+2\,z^3\right)
      }\over{2\,z^4}}&-{{7\,z^2-1}\over{4\,z^2}}&-{{1}\over{4}}\\
    {{
        \left(z-1\right)\,\left(1+z\right)\,\left(-1-3\,z-5\,z^2+z^3\right)
      }\over{2\,z^5}}&{{\left(z-1\right)\,\left(5\,z^2-1\right)}\over{4\,z
        ^3}}&{{z-1}\over{4\,z}}\\ 
    {{\left(z-1\right)^2\,\left(1+z\right)^4 }\over{2\,z^6}}&0&0 \\
  \end{bmatrix}
  $$
  which yields a $C^2$--convergent subdivision
  scheme that does not satisfy the spectral condition, but a
  generalized one with respect to the Taylor operator $\widetilde
  T_\Delta$. The result is shown in Fig.~\ref{fig:Exdelta}.
\end{example}

\noindent
For some time it was conjectured that all $C^d$ convergent Hermite
subdivision schemes must satisfy a spectral condition. This is
disproved by the following example of a family of convergent
schemes that satisfies no spectral condition.

\begin{theorem}
  If the nonzero elements of the matrix $\widetilde \Bb^*$ are of the
  form
  \begin{eqnarray*}
    \widetilde b_{jk}^* (z) & = & (z^{-1} - 1)^{j+1} \, h_{jk}^* (z), \qquad
                                  0 \le k < j < d, \\
    \widetilde b_{jj}^* (z) & = & \frac{(z^{-1} - 1)^{j+1}}{2^{j+1}},
                                  \qquad j=0,\dots,d-1, \\
    \widetilde b_{dj}^* (z) & = & \frac12 ( z^{-1} + 1 ) \, \left(
                                  z^{-2} - 1 \right)^{d-j} \qquad
                                  j=0,\dots,d, 
  \end{eqnarray*}
  then there exists a $C^d$--convergent Hermite subdivision scheme
  whose mask $\Ab$ satisfies $\widetilde T_\Delta S_\Ab =
	2^{-d}S_{\widetilde \Bb} \widetilde T_\Delta$. 
\end{theorem}

\begin{pf}
  Since $\Bb$ is lower triangular with contractions on the diagonal,
  the scheme $S_{\widetilde B}$ is contractive. The factorization is
  satisfied by construction.
\end{pf}

\subsection{A generic construction for arbitrary Taylor operators}
\label{sec:genconst}
For an arbitrary generalized Taylor operator $\widetilde T$, we
want to construct convergent schemes that factorize with respect to
$\widetilde T$, thus showing that convergence theory widely exceeds
spectral conditions.

\begin{theorem}\label{T:GenTaylorConvEx}
  For any $d \in \NN$ and any generalized Taylor operator $\widetilde
  T$ of order $d$ there exists a convergent Hermite subdivision scheme
  with mask $\Ab$ that \emph{factors} with $\widetilde T$, that is,
  such that $\widetilde T S_\Ab = 
	2^{-d}S_{\widetilde \Bb} 
  \widetilde T$ for some appropriate $\widetilde \Bb$.
\end{theorem}

\noindent
The proof continues the construction from the preceding subsection by
giving an explicit way to construct the polynomials $h_{dj}^*$,
$j=0,\dots,d$, in such a way that $S_\Ab$ admits the factorization.

\begin{pf}
  We will again set
  \begin{equation}
    \label{eq:bhconnect}
    \widetilde b_{dj}^* (z) = ( z^{-1} - 1 )^{d-j} \, h_{dj}^* (z^{-1}),
    \qquad 
  \end{equation}
  and make use of \eqref{eq:DelthCond} and \eqref{eq:DelthCond2} to
  determine the vectors 
  $$
  \hb_j =
  \begin{bmatrix}
    h_{j,j+1} \\ \vdots \\ h_{j1}
  \end{bmatrix}
  :=
  \begin{bmatrix}
    \left( h_{dj}^* \right)^{(j+1)} (1) \\ \vdots \\
    \left( h_{dj}^* \right)' (1)
  \end{bmatrix}
  \in \RR^{j+1}, \qquad j=0,\dots,d-1,
  $$
  which define $\widetilde \Bb^*$ and eventually
  the desired mask $\Ab^*$. We stack these vectors into
  the column vector
  $$
  \hb :=
  \begin{bmatrix}
    \hb_{d-1} \\ \vdots \\ \hb_0
  \end{bmatrix}
  \in \RR^{\frac{d(d+1)}{2}}.
  $$
  Again, let $h_{dd}^* (z)$ be the
  symbol of a contractive mask and recall that
  \begin{equation}
    \label{eq:TGenTaylorConvExPfMatType0}
    h_{dj}^* (1) = 2^{d-j}, \qquad j=0,\dots,d,
  \end{equation}
  is necessary due to Lemma~\ref{L:BtildReqs} to obtain $S_\Bb$ as a convergent
  vector subdivision scheme. Taking
  \eqref{eq:TGenTaylorConvExPfMatType0} into account, the requirement
  for $\hb_{d-1}$ can be obtained by setting $j=d$ in
  \eqref{eq:qDerivCondr}, which yields
  \begin{eqnarray*}
    \lefteqn{h_{d-1,r} + \sum_{k=1}^{r-1} w_{d,d-k} \frac{r!}{(r-k)!}
    h_{d-1-k,r-k}} \\
    & = & 2 \left( h_{dd}^* \right)^{(r)} (1) + r \left(
    h_{dd}^* \right)^{(r-1)} (1) - w_{d,d-r} 2^{r+1}, \qquad r=1,\dots,d-1.
  \end{eqnarray*}
  In the same way, \eqref{eq:qDerivCondj} transforms into
  $$
  h_{d-1,d} + \sum_{k=1}^{d-1} w_{d,d-k} \frac{d!}{(d-k)!}
    h_{d-1-k,d-k} = 2 \left( h_{dd}^* \right)^{(d)} (1) + d \, \left(
    h_{dd}^* \right)^{(d-1)} (1).
  $$
  In matrix form, this can be rewritten as
  \begin{eqnarray}
    \label{eq:TGenTaylorConvExPfMatType1}
    \bb_d
    & = & \left[
    \begin{array}{cccc|ccc|c|c}
      1 & & & & * & & & \dots & * \\
        & \ddots & & & & \ddots & & & 0 \\
        & & 1 & & & & * & & \vdots \\
        & & & 1 & 0 & \dots & 0 & \dots & 0
    \end{array}
                                          \right] \hb \\
    \nonumber
    & =: &
           \begin{bmatrix}
             \Ib_d & -\Hb_{d,d-2} & \dots & -\Hb_{d,0}
           \end{bmatrix} \hb,
  \end{eqnarray}
  where
  $$
  \Hb_{d,k} = -w_{d,k+1} \, 
  \begin{bmatrix}
    \frac{d!}{(k+1)!} \\
    & \ddots \\
    & & \frac{(d-k)!}{1!} \\
    0 & \dots & 0 \\
    \vdots & \ddots & \vdots \\
    0 & \dots & 0
  \end{bmatrix} \in \RR^{d \times k+1}, \qquad k = 0,\dots,d-2,
  $$
  and
  $$
  \bb_d :=
  \begin{bmatrix}
    2 \left( h_{dd}^* \right)^{(d)} (1) + d \, \left(
      h_{dd}^* \right)^{(d-1)} (1) \\
     2 \left( h_{dd}^* \right)^{(d-1)} (1) + (d-1) \left(
      h_{dd}^* \right)^{(d-2)} (1) - 2^d \, w_{d,1}  \\
    \vdots \\
    2 \left( h_{dd}^* \right)^{(1)} (1) +
    1 - 4 \, w_{d,d-1}
  \end{bmatrix} \in \RR^d.
  $$
  The conditions \eqref{eq:qDerivCondr} and \eqref{eq:qDerivCondj} for
  $q_{d-1}$ can, in the same way, be written as 
  $$
  0 = 2 h_{d-1,d-1} + (d-1) h_{d-1,d-2} - h_{d-2,d-1} -
  \sum_{k=1}^{d-2} w_{d-1,k} \frac{(d-1)!}{k!}
  h_{k-1,k},
  $$
  as well as for $r=2,\dots,d-2$,
  $$
  2^{r+2} \, w_{d-1,d-1-r} = 2 h_{d-1,r} + r \, h_{d-1,r-1} - h_{d-2,r} -
  \sum_{k=1}^{r-1} w_{d-1,d-1-k} \frac{r!}{(r-k)!}
  h_{d-2-k,r-k},
  $$
  and finally the case $r=1$
  $$
  2^{d-1} = 2 h_{d-1,r} - h_{d-2,r}.
  $$
  taking the matrix form
  \begin{eqnarray}
    \label{eq:TGenTaylorConvExPfMatType2}
    \bb_{d-1}
    & = & \left[
          \begin{array}{ccccc|cccc|c|c}
            0 & 2 & d-1 &&& -1 &&&& \dots & * \\
            0 & & 2 & \ddots &&& -1 &&&& 0 \\
            \vdots & & & \ddots & 2 &&& \ddots &&& \vdots \\
            0 & & & & 2 &&&& -1 & \dots & 0 \\
          \end{array}
    \right] \ab \\
    \nonumber
    & =: &
           \begin{bmatrix}
             \Cb_{d-1} & -\Ib_d & \Hb_{d-1,d-3} & \dots & \Hb_{d-1,0}
           \end{bmatrix} \hb
  \end{eqnarray}
  with
  $$
  \Cb_j :=
  \begin{bmatrix}
    0 & 2 & j  \\
    0 & & 2 & \ddots \\
    \vdots & & & \ddots & 2  \\
    0 & & & & 2  \\
  \end{bmatrix} \in \RR^j, \qquad j=2,\dots,d-1, \qquad
  \Cb_1 :=
  \begin{bmatrix}
    0 & 2 \\
    0 & 0 \\
  \end{bmatrix},
  $$
  and
  $$
  \Hb_{d-1,k} = -w_{d-1,k+1} \begin{bmatrix}
    \frac{(d-1)!}{(k+1)!} \\
    & \ddots \\
    & & \frac{(d-1-k)!}{1!} \\
    0 & \dots & 0 \\
    \vdots & \ddots & \vdots \\
    0 & \dots & 0
  \end{bmatrix} \in \RR^{d-1 \times k+1}, \qquad k = 0,\dots,d-3.
  $$
  With the general definition
  \begin{equation}
    \label{eq:TGenTaylorConvExPf1}
    \Hb_{j,k} = -w_{j,k+1} \begin{bmatrix}
    \frac{j!}{(k+1)!} \\
    & \ddots \\
    & & \frac{(j-k)!}{1!} \\
    0 & \dots & 0 \\
    \vdots & \ddots & \vdots \\
    0 & \dots & 0
  \end{bmatrix} \in \RR^{j \times k+1}, \qquad k = 0,\dots,j-2, \quad
  j=1,\dots,d,
  \end{equation}
  the conditions \eqref{eq:qDerivCondr} and \eqref{eq:qDerivCondj}
  result in the system
  \begin{equation}
    \label{eq:TGenTaylorConvExPf2}
    \bb = \begin{bmatrix}
      \bb_d \\ \vdots \\ \bb_1
    \end{bmatrix}
    = \begin{bmatrix}
      \Ib_d & -\Hb_{d,d-2} & -\Hb_{d,d-3} & \dots & -\Hb_{d,0} \\
      \Cb_{d-1} & -\Ib_d & \Hb_{d-1,d-3} & \dots & \Hb_{d-1,0}
      \\
      & \ddots & \ddots & \ddots & \vdots \\
      & & \Cb_1 & -1 & \Hb_{1,0}
    \end{bmatrix}
    \hb =: \Hb \hb.
  \end{equation}
  By Lemma~\ref{L:ADet}, which we prove next, this linear system has a
  unique solution $\hb$ for 
  any given polynomial $h_{dd}^* (z)$, which, by \eqref{eq:bhconnect}, defines
  the symbols $\widetilde b_{dj}^* (z)$, $j=0,\dots,d$, with
  $\widetilde b_{dd}^* (z) = h_{dd}^* (z)$ and therefore
  $$
  \widetilde \Bb^* (z) =
  \begin{bmatrix}
    \frac{z^{-1}-1}{2} \\
    (z^{-1}-1)^2 \, h_{10}^* (z) & \frac{(z^{-1}-1)^2}{4} \\
    \vdots & \ddots & \ddots \\
    (z^{-1}-1)^d \, h_{d-1,0}^* (z) & \dots & (z^{-1}-1)^d \,
    h_{d-1,d-2}^* (z) & \frac{(z^{-1}-1)^d}{2^d} \\
    \widetilde b_{d0}^* (z) & \dots & \widetilde b_{d,d-2}^* (z) &
    \widetilde b_{d,d-1}^* (z) & \widetilde b_{dd}^* (z)
  \end{bmatrix}
  $$
  is the symbol of a contractive scheme that satisfies the conditions
  from Lemma~\ref{L:BtildReqs} and for which there exists a mask $\Ab$
  such that $\widetilde S_\Ab = S_{\widetilde \Bb} \widetilde
  T$. Therefore, $\Ab$ defines a $C^d$--convergent Hermite subdivision
  scheme.
\end{pf}

\begin{remark}
  Recall that the whole construction process only had the purpose of
  finding the \emph{last row} of the lower triangular symbol
  $\widetilde B^*(z)$. All other entries could be chosen in a
  straightforward way.
\end{remark}

\begin{lemma}\label{L:ADet}
  Matrix $\Hb$ from \eqref{eq:TGenTaylorConvExPf2} satisfies $|\det
  \Hb | = 1$.
\end{lemma}

\begin{pf}
  Since the first column of $\Cb_{d-1}$ is zero, we can start with an
  expansion with respect to the first column, yielding that $\det \Hb$
  is the same as the determinant of $\Ab$ with first row and column
  erased. Then, we note that the last row of the matrix in
  \eqref{eq:TGenTaylorConvExPfMatType1} has only one nonzero entry,
  namely $-1$. Expansion with respect to this row also removes the
  column that contains the $2$ in the last row of
  \eqref{eq:TGenTaylorConvExPf2}. Expanding with respect to this row
  then removes the row that contains the last nonzero element in
  $\Hb_{d,d-2}$ in \eqref{eq:TGenTaylorConvExPfMatType1}, so that we
  can now expand with respect to the second last row of
  \eqref{eq:TGenTaylorConvExPfMatType1}. Circling in this way, we
  expand the determinant by means of factors that are $\pm 1$, hence,
  the determinant of $\Hb$ is $\pm 1$ and in particular independent of
  $\widetilde T$, that is, independent of $\wb_1,\dots,\wb_d$.
\end{pf}


\section{Examples}
\label{sec:Examples}

To illustrate the potential of the methods, we start with two examples
of masks obtained by the construction 
process in Theorem~\ref{T:GenTaylorConvEx}. We restrict ourselves to
the simplest nontrivial case $d=2$ here.

\begin{figure}
  \centering
  \includegraphics[width=.9\hsize]{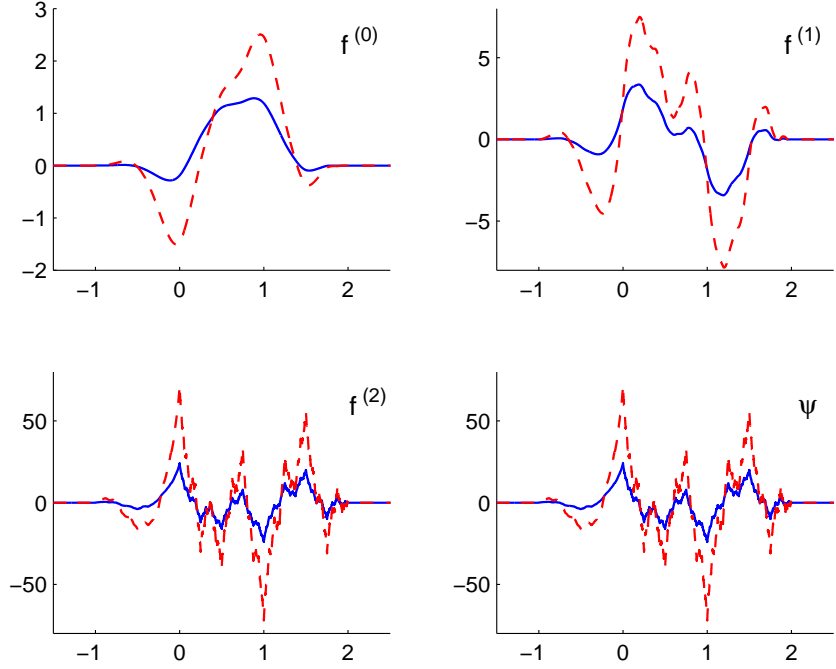}
  \caption{Limit functions for the constructions of
    Example~\ref{Ex:d=2arbw} for the values $w_{21} = \frac12$ (blue,
    solid) and $w_{21} = 1$ (red, dashed).}
  \label{fig:Exa1}
\end{figure}

\begin{example}\label{Ex:d=2arbw}
  One parameter, $w_{21}$,
  can be chosen freely. The associated linear system for $\hb$ becomes
  $$
  \begin{bmatrix}
    1 & & 2 w_{21} \\
    & 1 & \\
    & 2 & -1 
  \end{bmatrix}
  \begin{bmatrix}
    h_{12} \\ h_{11} \\ h_{01}
  \end{bmatrix}
  =
  \begin{bmatrix}
    2 \left( h_{dd}^* \right)'' (1) + 2 \left( h_{dd}^* \right)' (1)
    \\
    2 \left( h_{dd}^* \right)' (1) + 1 - 4 w_{21} \\
    2
  \end{bmatrix}
  $$
  which gives
  \begin{eqnarray*}
    h_{11} & = & 2 \left( h_{dd}^* \right)' (1) + 1 - 4 w_{21} \\
    h_{01} & = & 2 a_{11} - 2 = 4 \left( h_{dd}^* \right)' (1) - 8
                 w_{21} \\
    h_{12} & = & 2 \left( h_{dd}^* \right)'' (1) + 2 \left( h_{dd}^*
                 \right)' (1) - 2 w_{21} a_{01} \\
    & = & 2 \left( \left( h_{dd}^* \right)'' (1) + \left( h_{dd}^*
                 \right)' (1) \left( 1 - 4 w_{21} \right) + 8 w_{21}^2
          \right).
  \end{eqnarray*}
  Using the simplest possible choice $h_{dd}^* (z) = \frac12 ( z+1 )$,
  we get
  \begin{eqnarray*}
    h_{12} & = & 1 - 4 w_{21} + 16 w_{21}^2 \\
    h_{11} & = & 2 - 4 w_{21} \\
    h_{01} & = & 2 - 8 w_{21},
  \end{eqnarray*}
  and therefore
  \begin{eqnarray*}
    h_{21}^* (z)
    & = & \frac{\left( (1-4 w_{21}) z + (1+4 w_{21}) \right)^2}{2} + 2 w_{21} (
          z^2 - 1 ) \\
    h_{20}^* (z) & = & 4 + ( 2 - 8 w_{21} ) (z-1) = 2 \left( ( 1 - 4 w_{21} ) z
                       + (1 + 4 w_{21}) \right),
  \end{eqnarray*}
  yielding
  \begin{eqnarray*}
    \widetilde b_{22}^* (z) & = & \frac12 ( z^{-1} + 1 ) \\
    \widetilde b_{21}^* (z) 
                            & = & ( 1 - 4 w_{21} + 8 w_{21}^2 ) z^{-3}
                                  + 8 w_{21} ( 1 - 3 w_{21} ) z^{-2}
                                  - ( 1 + 4 w_{21} - 24 w_{21}^2 )
                                  z^{-1} - 8 w_{21}^2 \\
    \widetilde b_{20}^* (z) & = & ( 4 - 8 w_{21} ) z^{-2} - ( 4 - 16
                                  w_{21} ) z^{-1} + 8 w_{21}^2.
  \end{eqnarray*}
  The resulting limit functions are plotted in Fig~\ref{fig:Exa1}.
\end{example}

\begin{figure}
  \centering
  \includegraphics[width=.9\hsize]{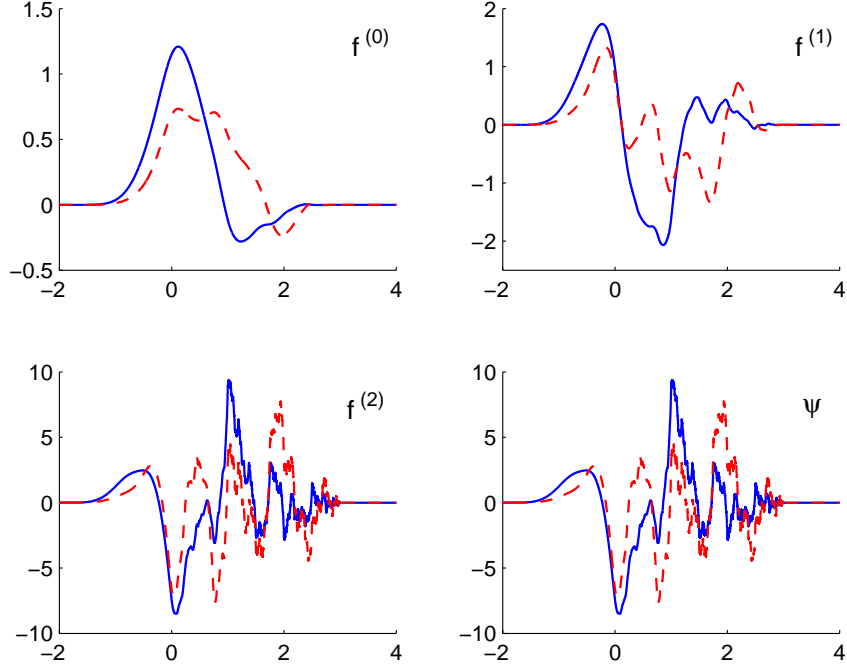}
  \caption{Limit functions for Example~\ref{Ex:d=2Bspl} for the values
    $w_{21} = \frac12$ (blue, solid) and $w_{21} = 1$ (red, dashed)
    and $n=5$.}
  \label{fig:Exan5}
\end{figure}

\begin{example}\label{Ex:d=2Bspl}
  In continuation of Example~\ref{Ex:d=2arbw}, we now choose an
  arbitrary contractive version based on
  $$
  h_{dd}^* (z) = \frac{(z+1)^n}{2^n}
  $$
  which has the property that
  $$
  h_{dd}^* (1) = 1, \qquad \left( h_{dd}^* \right)' (1) = \frac{n}2,
  \qquad \left( h_{dd}^* \right)'' (1) = \frac{n(n-1)}4,
  $$
  so that
  \begin{eqnarray*}
    h_{12} & = & 2 \left( \frac{n(n-1)}4 + \frac{n}2 \, ( 1-4 \,
                 w_{21}) + 8 w_{21}^2 \right)
                 = \frac{n(n+1)}2 - 4 n w_{21} + 16 w_{21}^2, \\
    h_{11} & = & n + 1 - 4 w_{21} \\
    h_{01} & = & 2n - 8 w_{21},
  \end{eqnarray*}
  which leads to the graphs shown in Fig.~\ref{fig:Exan5}. This even
  gives a whole family of convergent schemes with the additional
  parameter $n$.
\end{example}

\noindent
The last example revisits a Hermite subdivision scheme based on
B--splines that was introduced in \cite{MerrienSauer2011} and further
studied in \cite{MerrienSauer2017} as one of the first examples of a
family of convergent Hermite subdivision schemes that do not satisfy
the spectral condition.

This scheme is based on a construction detailed by Micchelli
in \cite{michelli95}.
Let 
$\varphi_0 (x) = \chi_{[0,1]}$
and define, for $r=1,2,\ldots$, 
the \emph{cardinal B--spline}
$\varphi_r=\varphi_0*\varphi_{r-1}$. 
We recall that  
$\varphi_r$ is a $C^{r-1}$ piecewise polynomial of degree $r$ 
with support $[0,r+1]$ that satisfies the refinement equation
$$
\varphi_r(x)=\frac1{2^r}\sum_{\alpha\in\ZZ}\begin{pmatrix}r+1\\ 
  \alpha\end{pmatrix}\varphi_r(2x-\alpha), \qquad
\begin{pmatrix}i\\ j\end{pmatrix}=\left\{ \begin{array}{ccl}\frac{i!}{j!(i-
j)!}
    & \; & \text{ if } 0\le j \le i,\\
    0 & & \text{otherwise}.
 \end{array}\right.
$$ 
The function
$v(x)=\sum_{\alpha\in\ZZ} f^{(0)}_0(\alpha) \varphi_r(x-\alpha)$
can be written as
$v(x)=
\sum_{\alpha\in\ZZ}f^{(0)}_n (\alpha) 
\varphi_r \left( 2^n \,x-\alpha \right)$, $n \in\NN_0$,  where 
\begin{equation}
  \label{subdivcoef0}
  f^{(0)}_{n+1} (\cdot) = 
  \sum_{\beta \in \ZZ}
  a_r(\cdot-2\beta)f^{(0)}_n(\beta),
  \qquad
  a_r (\alpha) = \frac{1}{2^r} { r+1 \choose \alpha }, \qquad \alpha
  \in \ZZ.
\end{equation}
We have proved in \cite[Proposition 5.3]{MerrienSauer2017} that  
for $i=0,\ldots,r$ one has
  \begin{equation}\label{eq:eigenpr}
    S_{a_r}p_i=\frac1{2^{i}}p_i, \qquad p_i:=\ell_r^{(r-i)}, \qquad
    \ell_r(x) := \frac{1}{r!}\prod_{j=1}^r(x+j).
\end{equation}
Taking derivatives of $v$,
$$\dfrac{d^iv}{dx^i}(x)=\sum_{\alpha\in\ZZ}2^{ni}
  \Delta^i f^{(0)}_n(\alpha-i) \, \varphi_{r-i} \left( 2^n x- \alpha
  \right),
  \qquad i = 0,\dots,r-1,$$ 
we define Hermite subdivision schemes of degree $d\le r$ with mask 
$\Ab(\alpha)$ and support $[0,r+d+1]$ by
applying differences to the mask $a_r$, yielding the following observation.

\begin{example}\label{Ex:SplineChain}
  The Hermite subdivision scheme based on
  $$
  \Ab (\alpha) =\begin{bmatrix}a_r(\alpha)&0&\ldots&0\\
    \Delta a_r(\alpha-1)&0&\ldots&0\\
    \Delta^2 a_r(\alpha-2)&0&\ldots&0\\
    \vdots\\ \Delta^da_r(\alpha-d)&0&\ldots&0\end{bmatrix}, \qquad
  \Ab^* (z) =\frac{(1+z)^{r+1}}{2^r}\begin{bmatrix}1&0&\ldots&0\\
    (1-z)&0&\ldots&0\\
    (1-z)^2&0&\ldots&0\\
    \vdots\\ (1-z)^d&0&\ldots&0\end{bmatrix}.
  $$
  has as limit function the vector consisting of the B--spline and its
  derivatives but does not satisfy the classical spectral condition,
  see \cite{MerrienSauer2011}.
\end{example}

\noindent
In the following, we prove that the Hermite scheme from
Example~\ref{Ex:SplineChain} possesses a spectral chain.

Firstly, the computation of Taylor expansions yields that
there for $p\in\Pi_d$ the vectors
$\vb_p=[p,p',\ldots,p^{(d)}]^T$ 
and $\hat \vb_p=[p,\Delta p(\cdot-1),\ldots,\Delta^d p(\cdot-d)]^T$
satisfy
$$
\hat \vb_p=\Rb \vb_p, \qquad \Rb :=
\begin{bmatrix}
  1&0&0&\ldots&0 \\
  &1& * & \dots & * \\
  & &\ddots& \ddots & \vdots \\ &&&1 & * \\
  & & & & 1
\end{bmatrix} \in\RR^{(d+1)\times(d+1)},
$$
where the $d-j$-th last components of $\hat \vb_p$ are zero if
$p\in\Pi_j$, $j<d$.  

Secondly, \eqref{eq:eigenpr} yields $S_{a_r}p_j=2^{-j}p_j$ and
the first component of $\vb_{p_j}$ is $p_j$ 
, since the only non zero column of the
matrices $\Ab (\alpha)$ is the first one, we therefore deduce that 
$$
S_{\Ab}\vb_{p_j} = S_{\Ab}
\begin{bmatrix} p_j\\ *
\end{bmatrix}
= S_{\Ab}\hat \vb_{p_j}=\frac1{2^j}\hat \vb_{p_j}, \qquad j=0,\dots,d,
$$
so that for $j=0,\ldots, d$, the vectors $\hat \vb_j = \hat \vb_{p_j}$
satisfy the spectral condition. To show that the associated 
$\hat \vb_j$ 
form a chain, we have to find the appropriate generalized Taylor
operator annihilating 
$\hat \vb_d$
, its uniqueness being guaranteed by
Corollary~\ref{C:DualTaylor}. This operator
is $\widetilde T_{S,d}$ from \eqref{eq:ChainExTSpline} in
Example~\ref{Ex:ChainEx}. Indeed, by 
Lemma~\ref{L:DiffIdent} 
proved 
at the end of this section,
\begin{eqnarray*}
  \left( \widetilde T_{S,d} \vb_d \right)_{d-j}
  & = & \Delta \left( \Delta^j p_d (\cdot-j)
        \right) - \sum_{k=1}^{d-j} \Delta^k \left( \Delta^j  p_d
        (\cdot - j) \right) (\cdot - k) \\
  & = & \Delta^d p_d ( \cdot - j - d + 1 ) - \Delta^d p_d ( \cdot -
        j -d ) = 0, \qquad j=0,\dots,d,
\end{eqnarray*}
since $\Delta^d p_d = 1$. The same argument also
shows that $\widetilde T_{S,d} \hat\vb_j = 0$, $j=0,\dots,d-1$.
Therefore $\Vb$
forms a spectral chain for $S_\Ab$ and by Theorem \ref{T:SpecChainFact} 
there  exists a finite mask
$\Bb \in \ell^{(d+1) \times (d+1)} (\ZZ)$ such that
$\widetilde T_{S,d} \, S_\Ab = S_{\widetilde \Bb} \, \widetilde T_{S,d}$.
  
\begin{example}[Example~\ref{Ex:SplineChain}
  continued]\label{Ex:SplineChain2}
  For $r=4$, $d=3$, we obtain
  $$
  \tilde \Bb^* (z)
  = \begin{bmatrix}
    -{{\left(z-1\right)^3\,z\,\left(1+z\right)^4}\over{2}}&{{
        \left(z-1\right)^2\,z^3\,\left(1+z\right)^3}\over{2}}&-{{\left(z-1
        \right)\,z^3\,\left(1+z\right)^2}\over{2}}&{{z^3\,\left(1+z\right)
      }\over{2}}\\[1mm]
    -{{\left(z-1\right)^3\,z\,\left(1+z\right)^4}\over{2}}
    &{{\left(z-1\right)^2\,z^3\,\left(1+z\right)^3}\over{2}}&-{{\left(z-
          1\right)\,z^3\,\left(1+z\right)^2}\over{2}}&{{z^3\,\left(1+z\right)
      }\over{2}}\\[1mm]
    -{{\left(z-1\right)^3\,z\,\left(1+z\right)^4}\over{2}}
    &{{\left(z-1\right)^2\,z^3\,\left(1+z\right)^3}\over{2}}&-{{\left(z-
          1\right)\,z^3\,\left(1+z\right)^2}\over{2}}&{{z^3\,\left(1+z\right)
      }\over{2}}\\[1mm]
    -{{\left(z-1\right)^3\,z\,\left(1+z\right)^4}\over{2}}
    &{{\left(z-1\right)^2\,z^3\,\left(1+z\right)^3}\over{2}}&-{{\left(z-
          1\right)\,z^3\,\left(1+z\right)^2}\over{2}}&{{z^3\,\left(1+z\right)
      }\over{2}}\\
  \end{bmatrix}.
  $$
\end{example}

\noindent
We close the paper with a simple identity on forward and backward
differences needed for Example~\ref{Ex:SplineChain2} that may,
however, be of independent interest.

\begin{lemma}\label{L:DiffIdent}
  For $p \in \Pi$ and $n \in \NN$ we have that
  \begin{equation}
    \label{eq:DiffIdent}
    \Delta p = \sum_{k=1}^{n-1} \Delta^k p (\cdot - k) + \Delta^n
    p(\cdot - n+1).
  \end{equation}
\end{lemma}

\begin{pf}
  Expanding the differences as
  $$
  \Delta^k p (\cdot - k) = \sum_{j=0}^k (-1)^j {k \choose j} p( \cdot
  - j ),
  $$
  we find that
  \begin{eqnarray*}
    \lefteqn{\Delta^n p(\cdot - n+1) + \sum_{k=1}^{n-1} \Delta^k p
    (\cdot - k) } \\ 
    & = &  p(\cdot + 1) + \sum_{j=0}^{n-1} (-1)^{j+1} { n \choose j+1}
          p(\cdot - j) + \sum_{k=1}^{n-1} \sum_{j=0}^k (-1)^j {k
          \choose j} p( \cdot - j ) \\ 
    & = & p( \cdot + 1) - p(\cdot) + \sum_{j=0}^{n-1} (-1)^j p( \cdot - j) \left(
          {n \choose j+1} - \sum_{k=j}^{n-1} {k \choose j} \right),
  \end{eqnarray*}
  from which the claim follows by taking into account the
  combinatorial identity
  \begin{equation}
    \label{eq:DiffIdentPf1}
    {n \choose j+1} = \sum_{k=j}^{n-1} {k \choose j}, \qquad 0 \le j
    \le n-1,
  \end{equation}
  which is easily proved by induction on $n$: calling the left
  hand side of \eqref{eq:DiffIdentPf1} $f(n)$ and the right hand side
  $g(n)$, the initial step $f(j+1) = g(j+1) = 1$ is obvious, while
  $$
  f(n+1) - f(n) = {n+1 \choose j+1} - {n \choose j+1} = {n \choose j}
  = g(n+1) - g(n)
  $$
  advances the induction.
\end{pf}



\end{document}